\documentclass[12pt,fleqn,leqno]{article}
\usepackage{amsmath,amssymb}

\usepackage[dvips]{graphicx}
\usepackage[dvips]{color}
\usepackage{lastpage}

\allowdisplaybreaks

\numberwithin{equation}{section}

\makeatletter
\renewcommand{\section}{\@startsection{section}{1}{0pt}{20pt}{6pt}{\large\bf}}
\renewcommand{\@seccntformat}[1]{\csname the#1\endcsname.\ }

\def\footnoterule{\kern -3pt \hrule width 2.7 true cm \kern 2.6pt}

\def\v{\vspace}
\def\h{\hspace}

\def\p{\!+\!}
\def\m{\!-\!}

\def\EE{\mathsf E\:\!}
\def\PP{\mathsf P}
\def\cF{{\cal F}}
\def\R{I\!\!R}
\def\LL{I\!\!L}
\def\eps{\varepsilon}

\begin{document}

\title{\bf Quickest Real-Time Detection of \\ Multiple Brownian Drifts}
\author{P. A. Ernst\footnote{Department of Statistics, Rice University. Email: philip.ernst@rice.edu}, H. Mei\footnote{Department of Statistics, Rice University. Email: hongwei.mei@rice.edu} \:\! \&\:\! G. Peskir\footnote{Department of Mathematics, The University of Manchester. Email: goran@maths.man.ac.uk}}
\date{}
\maketitle




{\par \leftskip=2cm \rightskip=2cm \footnotesize

Consider the motion of a Brownian particle in $n$ dimensions, whose
coordinate processes are standard Brownian motions with zero drift
initially, and then at some random/unobservable time, exactly $k$ of
the coordinate processes get a (known) non-zero drift permanently.
Given that the position of the Brownian particle is being observed
in real time, the problem is to detect the time at which the $k$
coordinate processes get the drift as accurately as possible. We
solve this problem in the most uncertain scenario when the
random/unobservable time is (i) exponentially distributed and (ii)
independent from the initial motion without drift. The solution is
expressed in terms of a stopping time that minimises the probability
of a false early detection and the expected delay of a missed late
detection. The elliptic case $k=1$ has been settled in \cite{EP}
where the hypoelliptic case $1 < k < n$ resolved in the present
paper was left open (the case $k = n$ reduces to the classic case
$n=1$ having a known solution). We also show that the methodology
developed solves the problem in the general case where exactly $k$
is relaxed to \emph{any} number of the coordinate processes getting
the drift. To our knowledge this is the first time that such a
multi-dimensional hypoelliptic problem has been solved exactly in
the literature.

\par}





\section{Introduction}

Imagine the motion of a Brownian particle in $n$ dimensions, whose
coordinate processes are standard Brownian motions with zero drift
initially, and then at some random/unobservable time $\theta$,
exactly $k$ of the coordinate processes get a (known) non-zero drift
$\mu$ permanently. Assuming that the position of the Brownian
particle is being observed in real time, the problem is to detect
the time $\theta$ at which the $k$ coordinate processes get the
drift $\mu$ as accurately as possible. In the most uncertain
scenario, where $\theta$ is assumed to be (i) exponentially
distributed and (ii) independent from the initial motion without
drift, the solution to the problem when $k=1$ has been derived in
\cite{EP}. The purpose of the present paper is to derive the
solution to the problem when $1 < k < n$ that was left open in
\cite{EP}. (Note that the case $k = n$ reduces to the classic case
$n=1$ having a known solution and is therefore excluded throughout.)

Denoting the position of the Brownian particle in $n$ dimensions by
$X$, the error to be minimised over all stopping times $\tau$ of $X$
is expressed as the the linear combination of the probability of the
\emph{false alarm} $\PP_{\!\pi}(\tau\! <\! \theta)$ and the expected
\emph{detection delay} $\EE_\pi(\tau \m \theta)^+$ where $\pi \in
[0,1]$ denotes the probability that $\theta$ has already occurred at
time $0$. This problem formulation of quickest detection dates back
to \cite{Sh-1} and has been extensively studied to date (see
\cite{Sh-4} and the references therein). The linear combination
represents the Lagrangian and once the optimal stopping problem has
been solved in that form this will also lead to the solution of the
constrained problems where an upper bound is imposed on either the
probability of the false alarm or the expected detection delay
respectively.

Quickest detection problems related to the present problem have been
studied by a number of authors and we refer to \cite[Section 1]{EP}
for an overview of this literature (see also Remark 9 in \cite{EP}).
The initial optimal stopping problem in these quickest detection
problems is equivalent to an optimal stopping problem for the
posterior probability distribution ratio process $\varPhi$ of
$\theta$ given $X$. The infinitesimal generator of the
Markov/diffusion process $\varPhi$ (combined with $X$ when needed)
in these optimal stopping problems is of \emph{parabolic} type. In
contrast, this was no longer the case in the present problem when
$k=1$ as noted in \cite{EP}, where the infinitesimal generator
$\LL_\varPhi$ of the multi-dimensional Markov/diffusion process
$\varPhi$ is of \emph{elliptic} type.

The present quickest detection problem when $1< k < n$ was left open
in \cite{EP}. We will see below that the ellipticity of
$\LL_\varPhi$ remains valid for $k = n \m 1$ but breaks down when $1
< k < n \m 1$. We show however that $\LL_\varPhi$ satisfies the
\emph{H\"ormander condition} (cf.\ \cite{Ho}) in this case so that
$\LL_\varPhi$ is \emph{hypoelliptic}. This fact may appear somewhat
surprising at first glance given that $\varPhi$ is governed by a
system of $N = \binom{n}{k}$ stochastic differential equations
driven by $n$ Brownian motions. For example, when $n=10$ and $k=5$
then $N=252$ so that the system governing $\varPhi$ consists of
$252$ stochastic differential equations driven by $10$ Brownian
motions. Nonetheless, we establish that the structure of this system
when $1< k < n \m 1$ is sufficiently supported to ensure
\emph{hypoellipticity} as a substitute for the broken ellipticity of
$k$ being either $1$ or $n \m 1$. This provides regularity of the
value function in the optimal stopping set. Moreover, the previous
conclusions on the space operator $\LL_\varPhi$ extend to the
backward time-space operator $-\partial_t \p \LL_\varPhi$ as well
which in turn imply that $\varPhi$ is a strong Feller process. This
opens avenues to regularity of the value function at the optimal
stopping points and closes our exploratory analysis of the problem.

In the reminder we show that the hypoelliptic structure of the
infinitesimal generator combined with the concavity of the loss
functional in the optimal stopping problem is sufficiently robust to
yield the solution. Finding the \emph{exact} solution to the
quickest detection problem for the observed process $X$ in $n$
dimensions when $1 < k < n$ is the main contribution of the present
paper. We also show that the methodology developed solves the
problem in the general case where exactly $k$ is relaxed to
\emph{any} number of the coordinate processes getting the drift. To
our knowledge this is the first time that such a multi-dimensional
hypoelliptic problem has been solved exactly in the literature.

\section{Formulation of the problem}

In this section we formulate the quickest detection problem under
consideration. The initial formulation of the problem will be
revaluated under a change of measure in the next section.

\medskip

1.\ We consider a Bayesian formulation of the problem where it is
assumed that one observes a sample path of the standard\!
$n$\!-dimensional Brownian motion $X=(X^1, \ldots, X^n)$, whose
coordinate processes $X^1, \ldots, X^n$ are standard Brownian
motions with zero drift initially, and then at some
random/unobservable time $\theta$ taking value $0$ with probability
$\pi \in [0,1]$ and being exponentially distributed with parameter
$\lambda>0$ given that $\theta>0$, exactly $k$ of the coordinate
processes $X^1, \ldots, X^n$ get a (known) non-zero drift $\mu$
permanently. The problem is to detect the time $\theta$ at which the
$k$ coordinate processes get the drift $\mu$ as accurately as
possible (neither too early nor too late). This problem belongs to
the class of quickest real-time detection problems as discussed in
Section 1 above.

\smallskip

2.\ The observed process $X=(X^1, \ldots, X^n)$ satisfies the
stochastic differential equations
\begin{equation} \h{5pc} \label{2.1}
dX_t^i = \mu\:\! I(i\! \in\! \beta,t\! \ge\! \theta)\;\! dt +
dB_t^i\;\; (1\! \le\! i\! \le\! n)
\end{equation}
driven by a standard\! $n$\!-dimensional Brownian motion $B=(B^1,
\ldots, B^n)$ under the probability measure $\PP_{\!\pi}$ specified
below, where the random variable $\beta$ taking values in the set
$C_k^n := \{ (n_1, \ldots, n_k)\, \vert\, 1 \le n_1 < \ldots < n_k
\le n \}$ satisfies $\PP_{\!\pi}(\beta\! =\! (n_1, \ldots,
n_k))=p_{n_1, \ldots, n_k}$ for some $p_{n_1, \ldots, n_k} \in
[0,1]$ with $\sum p_{n_1, \ldots, n_k} = 1$ given and fixed where
the sum is taken over all $(n_1, \ldots, n_k)$ in $C_k^n$. With a
slight abuse of notation, in \eqref{2.1} we write $i \in \beta$ to
express the fact that $i$ belongs to the set $\{ n_1, \ldots, n_k
\}$ consisting of the elements which form $\beta = (n_1, \ldots,
n_k)$ in $C_k^n$. This means that $n_1, \ldots, n_k \in \beta$ if
and only if the coordinate processes $X^{n_1}, \ldots, X^{n_k}$ get
drift $\mu$ at time $\theta$ with probability $p_{n_1, \ldots, n_k}$
for $(n_1, \ldots, n_k) \in C_k^n$. With a similar abuse of
notation, which will be helpful in what follows, we will arrange the
elements of $C_k^n$ in a lexicographic order and identify the\!
$i$\!-th element of the ordered set $C_k^n$ by its index $i$ itself
for $1 \le i \le N$ where we set $N:=\binom{n}{k}$ to denote the
total number of elements in $C_k^n$. Often we will write $i = (n_1,
\ldots, n_k)$ or $i \in C_k^n$ to express this identification
explicitly for $1 \le i \le N$ while at other places it will be
clear from the context whether the index $i$ belongs to the set
$C_k^n$ in this sense or the set $\{1, \ldots, n\}$ as in
\eqref{2.1} above. The unobservable time $\theta$, the unknown
coordinates $\beta$, and the driving Brownian motion $B$ are all
assumed to be independent under $\PP_{\!\pi}$ for $\pi \in [0,1]$
given and fixed.

\medskip
3.\ Standard arguments imply that the previous setting can be
realised on a probability space $(\Omega,\cF,\PP_{\!\pi})$ with the
probability measure $\PP_{\!\pi}$ being decomposable as follows
\begin{equation} \h{5pc} \label{2.2}
\PP_{\!\pi} = \sum_{i \in C_k^n} p_i \Big(\pi\;\! \PP_{\!i}^0
+ (1 \m \pi)\! \int_0^\infty\! \lambda\;\! e^{-\lambda t}\,
\PP_{\!i}^t\, dt \Big)
\end{equation}
for $\pi \in [0,1]$ where $\PP_{\!i}^t$ is the probability measure
under which the coordinate processes $X^{n_1}, \ldots, X^{n_k}$ of
the observed process $X$ get drift $\mu$ at time $t \in [0,\infty)$
for $i = (n_1, \ldots, n_k)$ in $C_k^n$ and $\lambda>0$ is given and
fixed. The decomposition \eqref{2.3} expresses the fact that the
unobservable time $\theta$ is a non-negative random variable
satisfying $\PP_{\!\pi}(\theta = 0) = \pi$ and $\PP_{\!\pi}(\theta >
t\, \vert\, \theta > 0) = e^{-\lambda t}$ for $t>0$. Thus
$\PP_{\!i}^t(X \in\, \cdot\, ) = \PP_{\!\pi}(X \in\, \cdot\; \vert\,
\beta=i,\:\! \theta=t)$ is the probability law of the standard\!
$n$\!-dimensional Brownian motion process $X$ whose coordinate
processes $X^{n_1}, \ldots, X^{n_k}$ get drift $\mu$ at time $t \in
[0,\infty)$ for $i = (n_1, \ldots, n_k)$ in $C_k^n$. To remain
consistent with this notation we also denote by $\PP_{\!i}^\infty$
the probability measure under which the coordinate processes
$X^{n_1}, \ldots, X^{n_k}$ of $X$ get no drift $\mu$ at a finite
time for $i = (n_1, \ldots, n_k)$ in $C_k^n$. Thus
$\PP_{\!i}^\infty(X \in\, \cdot\; ) = \PP_{\!\pi}(X \in\, \cdot\;
\vert\, \beta=i,\:\! \theta = \infty)$ is the probability law of the
standard\! $n$\!-dimensional Brownian motion process for $i \in
C_k^n$. Clearly the subscript $i$ is superfluous in this case and we
will often write $\PP^\infty$ instead of $\PP_{\!i}^\infty$ for $i
\in C_k^n$. Moreover, by $\PP_{\!i}$ we denote the probability
measure under which the coordinate processes $X^{n_1}, \ldots,
X^{n_k}$ of $X$ get drift $\mu$ at time $\theta$ for $i = (n_1,
\ldots, n_k)$ in $C_k^n$. From \eqref{2.2} we see that
\begin{equation} \h{8pc} \label{2.3}
\PP_{\!\pi} = \sum_{i \in C_k^n} p_i\:\! \PP_{\!i}
\end{equation}
where $\PP_{\!i} = \pi\;\! \PP_{\!i}^0 + (1 \m \pi)\!
\int_0^\infty\! \lambda\; \! e^{-\lambda t}\, \PP_{\!i}^t\, dt$ for
$i \in C_k^n$ and $\pi \in [0,1]$. Note that $\PP_{\!i}$ depends on
$\pi \in [0,1]$ as well but we will omit this dependence from its
notation for $i \in C_k^n$.

\medskip

4.\ Being based upon continuous observation of $X=(X^1, \ldots,
X^n)$, the problem is to find a stopping time $\tau_*$ of $X$ (i.e.\
a stopping time with respect to the natural filtration $\cF_t^X =
\sigma(X_s\, \vert\, 0 \le s \le t)$ of $X$ for $t \ge 0$) that is
`as close as possible' to the unknown time $\theta$. More precisely,
the problem consists of computing the value function
\begin{equation} \h{4pc} \label{2.4}
V(\pi) = \inf_\tau \Big[ \PP_{\!\pi}(\tau < \theta) + c\;\! \EE_\pi
(\tau - \theta)^+ \Big]
\end{equation}
and finding the optimal stopping time $\tau_*$ at which the infimum
in \eqref{2.4} is attained for $\pi \in [0,1]$ and $c>0$ given and
fixed \v{-1pt} (recalling also that $p_i \in [0,1]$ with $\sum_{i
\in C_k^n} p_i = 1$ are given and fixed). Note in \eqref{2.4} that
$\PP_{\!\pi}(\tau < \theta)$ is the probability of the \emph{false
alarm} and $\EE_\pi(\tau - \theta)^+$ is the expected
\emph{detection delay} associated with a stopping time $\tau$ of $X$
for $\pi \in [0,1]$. Recall also that the expression on the
right-hand side of \eqref{2.4} is the Lagrangian associated with the
constrained problems as discussed in Section 1 above.

\medskip

5.\ To tackle the optimal stopping problem \eqref{2.4} we consider
the \emph{posterior probability distribution process}
$\varPi=(\varPi_t)_{t \ge 0}$ of $\theta$ given $X$ that is defined
by
\begin{equation} \h{7pc} \label{2.5}
\varPi_t = \PP_{\!\pi}(\theta \le t\, \vert\, \cF_t^X)
\end{equation}
for $t \ge 0$. Note that we have
\begin{equation} \h{8pc} \label{2.6}
\varPi_t = \sum_{i=1}^N \varPi_t^i
\end{equation}
where the summands are defined as follows
\begin{equation} \h{6pc} \label{2.7}
\varPi_t^i = \PP_{\!\pi}(\beta=i,\: \theta \le t\, \vert\, \cF_t^X)
\end{equation}
for $t \ge 0$ and $1 \le i \le N$ in $C_k^n$. The right-hand side of
\eqref{2.4} can be rewritten to read
\begin{equation} \h{5pc} \label{2.8}
V(\pi) = \inf_\tau\;\! \EE_\pi \Big( 1 \m \varPi_\tau + c \int_0^\tau\!
\varPi_t\, dt \Big)
\end{equation}
for $\pi \in [0,1]$.

\medskip

6.\ To connect the process $\varPi$ to the observed process $X$ we
define the \emph{posterior probability distribution ratio process}
$\varPhi=(\varPhi^1, \ldots, \varPhi^N)$ of $\theta$ given $X$ by
\begin{equation} \h{8pc} \label{2.9}
\varPhi_t^i = \frac{\varPi_t^i}{\bar \varPi_t^i}
\end{equation}
where the denominator is defined is follows
\begin{equation} \h{7pc} \label{2.10}
\bar \varPi_t^i = \PP_{\!\pi}(\beta=i,\: \theta > t\, \vert\, \cF_t^X)
\end{equation}
for $t \ge 0$ and $1 \le i \le N$ in $C_k^n$. By the Girsanov
theorem we see that the \emph{likelihood ratio process} $L=(L^1,
\ldots L^N)$ can be expressed as follows
\begin{equation} \h{5pc} \label{2.11}
L_t^i = \frac{d \PP_{\!i,t}^0}{d \PP_{\!t}^\infty} = \exp \Big(
\mu \sum_{j=1}^k X_t^{n_j} - k\, \frac{\mu^2}{2}\:\! t \Big)
\end{equation}
for $t \ge 0$ and $1 \le i \le N$ identified with $(n_1, \ldots,
n_k)$ in $C_k^n$ where $\PP_{\!i,t}^0$ and $\PP_{\!t}^\infty$ denote
the restrictions of the measures $\PP_{\!i}^0$ and $\PP^\infty$ to
$\cF_t^X$ for $t \ge 0$ and $1 \le i \le N$ in $C_k^n$. Moreover,
using the same arguments as in \cite[Section 2]{EP} we find that
\begin{equation} \h{7pc} \label{2.12}
\varPhi_t^i = e^{\lambda t} L_t^i\:\! \Big( \varPhi_0^i + \lambda
\int_0^t \frac{ds} {e^{\lambda s} L_s^i}\, \Big)
\end{equation}
with $\varPhi_0^i=\pi/(1 \m \pi)$ for $t \ge 0$ and $1 \le i \le N$
in $C_k^n$. From \eqref{2.11} and \eqref{2.12} we see that the
process $\varPhi=(\varPhi^1, \ldots\, \varPhi^N)$ is an explicit
(path-dependent) functional of the observed process $X=(X^1, \ldots,
X^n)$ and hence observable (by observing a sample path of $X$ we are
also seeing a sample path of $\varPhi$ both in real time).

\section{Measure change}

In this section we show that changing the probability measure
$\PP_{\!\pi}$ for $\pi \in [0,1]$ to $\PP^\infty$ in the optimal
stopping problem \eqref{2.4} or \eqref{2.8} provides crucial
simplifications of the setting which make the subsequent analysis
possible. This will be achieved by invoking the decomposition of
$\PP_{\!\pi}$ into $\PP_{\!i}$ for $i \in C_k^n$ as stated in
\eqref{2.2} above, changing each probability measure $\PP_{\!i}$ to
$\PP_{\!i}^\infty$, and recalling that each $\PP_{\!i}^\infty$
coincides with $\PP^\infty$ for $i \in C_k^n$.

\medskip

\noindent 1.\ We show that the optimal stopping problem \eqref{2.8} admits a
transparent reformulation under the probability measure $\PP^\infty$
in terms of the process $\varPhi=(\varPhi^1, \ldots, \varPhi^N)$
defined in \eqref{2.9} above. Recall that $\varPhi^i$ starts at
$\pi/(1 \m \pi)$ and this dependence on the initial point will be
indicated by a superscript to $\varPhi^i$ when needed for $1 \le i
\le N$ in $C_k^n$.

\medskip

\textbf{Proposition 1.} \emph{The value function $V$ from
\eqref{2.8} satisfies the identity
\begin{equation} \h{7pc} \label{3.1}
V(\pi) = (1 \m \pi)\, \big[ 1 + c\;\! \hat V(\pi) \big]
\end{equation}
where the value function $\hat V$ is given by
\begin{equation} \h{3pc} \label{3.2}
\hat V(\pi) = \inf_\tau\;\! \EE^\infty \Big[ \int_0^\tau e^{-\lambda
t} \Big( \sum_{i=1}^N p_i\;\! \varPhi_t^{i,\pi/(1-\pi)} -
\frac{\lambda} {c}\;\! \Big)\;\! dt\:\! \Big]
\end{equation}
for $\pi \in [0,1)$ and the infimum in \eqref{3.2} is taken over all
stopping times $\tau$ of $X$.}

\smallskip
\textbf{Proof.} A derivation  of \eqref{3.1} and \eqref{3.2} can be
reduced to one dimension where the change-of-measure identity (4.12)
from \cite{JP} is applicable in exactly the same way as in the proof
of Proposition 1 in \cite{EP}. Given that the present
multi-dimensional setting creates no additional difficulties we will
omit details and this completes the proof. \hfill $\square$

\medskip

2.\ From Proposition 1 we see that the optimal stopping problem
\eqref{2.4} or \eqref{2.8} is equivalent to the optimal stopping
problem \eqref{3.2}. From \eqref{2.11}+\eqref{2.12} using It\^o's
formula we find that
\begin{equation} \h{3pc} \label{3.3}
d \varPhi_t^i = \lambda\:\! (1 \p \varPhi_t^i)\;\! dt + \sum_{j=1}^k
\mu\;\! \varPhi_t^i\, dB_t^{n_j}\;\; (i\! =\! (n_1, \ldots, n_k)\!
\in\! C_k^n)
\end{equation}
under $\PP^\infty$ with $\varPhi_0^i = \varphi_i$ in $[0,\infty)$
all being equal to $\pi/(1 \m \pi)$ for $1 \le i \le N$ in $C_k^n$
and $\pi \in [0,1)$. Moreover, the system of stochastic differential
equations \eqref{3.3} has a unique strong solution given by
\eqref{2.11}+\eqref{2.12} above, where $X$ equals $B$ under
$\PP^\infty$, so that $\varPhi = (\varPhi^1, \ldots, \varPhi^N)$ is
a \emph{strong Markov} process (see e.g.\ \cite[pp 158-163]{RW}). We
will establish further in Section 5 below that $\varPhi$ is a strong
Feller process but these arguments are more subtle at this stage.
Noticing that $W_t^i := (1/\sqrt{k}) \sum_{j=1}^k B_t^{n_j}$ defines
a standard Brownian motion process for $t \ge 0$, we see from
\eqref{3.3} that each $\varPhi^i$ is a Shiryaev diffusion process
for $1 \le i \le N$ identified with $(n_1, \ldots, n_k)$ in $C_k^n$.
Basic properties of the Shiryaev diffusion processes are reviewed in
\cite[Section 2]{Pe-2}. In particular, it is known that $\varPhi^i$
is \emph{recurrent} in $[0,\infty)$ if and only if $\lambda \le k
\mu^2\!/2$ for $1 \le i \le N$ in $C_k^n$. If $\lambda
> k \mu^2\!/2$ then $\varPhi^i$ is \emph{transient} in $[0,\infty)$
with $\varPhi_t^i \rightarrow \infty$ almost surely under
$\PP^\infty$ as $t \rightarrow \infty$ for $1 \le i \le N$ in
$C_k^n$.

\medskip

3.\ To tackle the equivalent optimal stopping problem \eqref{3.2}
for the strong Markov process $\varPhi = (\varPhi^1, \ldots,
\varPhi^N)$ solving \eqref{3.3} we will enable $\varPhi =
(\varPhi^1, \ldots, \varPhi^N)$ to start at any point $\varphi =
(\varphi_1, \ldots, \varphi_N) \in [0,\infty)^N$ under the
probability measure $\PP_{\!\varphi}^\infty$ so that the optimal
stopping problem \eqref{3.2} extends as follows
\begin{equation} \h{4.5pc} \label{3.4}
\hat V(\varphi) = \inf_\tau\;\! \EE_\varphi^\infty \Big[ \int
_0^\tau e^{-\lambda t} \Big( \sum_{i=1}^N p_i\;\! \varPhi_t^i
- \frac{\lambda} {c}\;\! \Big)\;\! dt\:\! \Big]
\end{equation}
for $\varphi \in [0,\infty)^N$ with $\PP_{\!\varphi}^\infty
(\varPhi_0\! =\! \varphi)=1$ \v{-1pt} where the infimum is taken
over all stopping times $\tau$ of $\varPhi$ and we recall that $p_i
\in [0,1]$ for $1 \le i \le N$ in $C_k^n$ with $\sum_{i=1}^N p_i =
1$. In this way we have reduced the initial quickest detection
problem \eqref{2.4} or \eqref{2.8} to the optimal stopping problem
\eqref{3.4} for the strong Markov process $\varPhi = (\varPhi^1,
 \ldots, \varPhi^N)$ solving \eqref{3.3} and being explicitly
given by the Markovian flow \eqref{2.11}+\eqref{2.12} of the initial
point $(\varPhi_0^1, \ldots, \varPhi_0^N) = (\varphi_1, \ldots
\varphi_N) =: \varphi$ in $[0,\infty)^N$ under
$\PP_{\!\varphi}^\infty$. Note that the optimal stopping problem
\eqref{3.4} is inherently/fully\!\! $N$\!-dimensional and the
infinitesimal generator of $\varPhi = (\varPhi^1, \ldots,
\varPhi^N)$ is of elliptic type when $k=1$ or $k = n \m 1$ (because
$N=n$ and the diffusion matrix in \eqref{3.3} is regular in either
case) but is not of elliptic (or parabolic) type when $1 < k < n \m
1$ (because $N>n$ and hence the diffusion matrix in \eqref{3.3}
cannot be regular). We will return to this issue in Section 5 below.

\section{Mayer formulation}

The optimal stopping problem \eqref{3.4} is Lagrange formulated. In
this section we derive its Mayer reformulation which is helpful in
the subsequent analysis.

1.\ From \eqref{3.3} we read that the infinitesimal generator of the
strong Markov process $\varPhi = (\varPhi^1, \ldots, \varPhi^N)$ is
given by
\begin{equation} \h{3pc} \label{4.1}
\LL_\varPhi = \sum_{i=1}^N \lambda\:\! (1 \p \varphi_i)\;\!
\partial_{\varphi_i} + \frac{1}{2} \sum_{i,j=1}^N \mu^2\:\!
\varphi_i \varphi_j\:\! (I_i,I_j)\:\! \partial_{\varphi_i \varphi_j}
\end{equation}
for $(\varphi_1, \ldots \varphi_N) \in (0,\infty)^N$ where $I_i =
(I_{i1}, \ldots, I_{in})$ with $I_{ip} = 1_i(p)$ which by definition
equals $1$ if $p$ belongs to $\{n_1, \ldots, n_k\}$ for $1 \le i \le
N$ identified with $(n_1, \ldots, n_k)$ in $C_k^n$ and $0$
otherwise, and similarly $I_j = (I_{j1}, \ldots, I_{jn})$ with
$I_{jp} = 1_j(p)$ which by definition equals $1$ if $p$ belongs to
$\{m_1, \ldots, m_k\}$ for $1 \le j \le N$ identified with $(m_1,
\ldots, m_k)$ in $C_k^n$ and $0$ otherwise, while $(I_i,I_j)$
denotes the scalar product of $I_i$ and $I_j$ given by $\sum_{p=1}^n
I_{ip} I_{jp}$ for $1 \le i,j \le N$ in $C_k^n$. From \eqref{2.12}
we see that the topological boundary $\{\, (\varphi_1, \ldots,
\varphi_N) \in [0,\infty)^N\; \vert\; \varphi_i=0\;\, \text{for
some}\;\, 1 \le i \le N\, \}$ of the state space $[0,\infty)^N$
consists of \emph{entrance} boundary points for $\varPhi$ (meaning
that $\varPhi$ can be started at any boundary point never to return
to the boundary) and clearly the differential operator $\LL_\varPhi$
is of elliptic type when $k=1$ or $k = n \m 1$ but is not of
elliptic (or parabolic) type when $1 < k < n \m 1$ which is the case
of main interest in what follows. We will determine its type when $1
< k < n \m 1$ in Section 5 below.\\

2.\ For the Mayer reformulation of the problem \eqref{3.4} we need
to look for a function $M : [0,\infty)^N \rightarrow \R$ solving the
partial differential equation
\begin{equation} \h{8pc} \label{4.2}
\LL_\varPhi M \m \lambda M = L
\end{equation}
on $(0,\infty)^N$ where in view of \eqref{3.4} we set
\begin{equation} \h{5pc} \label{4.3}
L(\varphi_1, \ldots, \varphi_N) = \sum_{i=1}^N p_i \varphi_i
- \lambda/c
\end{equation}
for $(\varphi_1, \ldots, \varphi_N) \in [0,\infty)^N$. Noting that
the mixed derivatives in \eqref{4.1} vanish, and ignoring the
constant $-\lambda/c$ on the right-hand side of \eqref{4.2} with
\eqref{4.3} for now, we see that a possible attempt to solve the
resulting partial differential equation is to separate the variables
$\varphi_i$ by considering the ordinary differential equations
\begin{equation} \h{5pc} \label{4.4}
\lambda\:\! (1 \p \varphi_i)\:\! M_i' + \frac{\mu^2}{2}\:\!
\varphi_i^2\: M_i'' - \lambda\:\! M_i = \varphi_i
\end{equation}
where $M_i = M_i(\varphi_i)$ is a function/solution to be found for
$\varphi_i \in (0,\infty)$ with $1 \le i \le N$ in $C_k^n$. Note
that $M_i = M_j$ for $i \ne j$ in $C_k^n$ so that the subscript to
$M$ is superfluous but we will keep it to distinguish $M_i$ for $1
\le i \le N$ in $C_k^n$ from the general function $M$ solving
\eqref{4.2} to be defined shortly below. It was shown in
\cite[Section 4]{EP} that the sought solution to \eqref{4.4} is
given by
\begin{equation} \h{2pc} \label{4.5}
M_i(\varphi_i) = \frac{2}{\mu^2}\, (1 \p \varphi_i)\! \int_0^{\varphi_i
/(1 + \varphi_i)}\!\! \Big( \frac{1 \m v}{v} \Big)^{\!\kappa} e^{\kappa
/v}\! \int_0^v \frac{u^{\kappa-1}}{(1 \m u)^{\kappa+2}}\, e^{-\kappa/u}
\, du\;\! dv
\end{equation}
for $\varphi_i \in [0,\infty)$ and $1 \le i \le N$ in $C_k^n$ where
we set $\kappa = 2 \lambda/\mu^2$. Define a function $M :
[0,\infty)^N \rightarrow \R$ by setting
\begin{equation} \h{5pc} \label{4.6}
M(\varphi_1, \ldots, \varphi_N) = \sum_{i=1}^N p_i\:\! M_i(
\varphi_i) + 1/c
\end{equation}
for $(\varphi_1, \ldots, \varphi_N) \in [0,\infty)^N$. The arguments
above then show that the function $M$ from \eqref{4.6} solves the
equation \eqref{4.2} above (notice that the final term $1/c$ yields
the missing constant $-\lambda/c$ on the right-hand side of
\eqref{4.2} with \eqref{4.3} as needed).

\medskip

3.\ Having defined the function $M$ in \eqref{4.6} we can now
describe the Mayer reformulation of the optimal stopping problem
\eqref{3.4} as follows.

\medskip

\textbf{Proposition 2.} \emph{The value function $\hat V$ from
\eqref{3.4} can be expressed as
\begin{equation} \h{5pc} \label{4.7}
\hat V(\varphi) = \inf_\tau\:\! \EE_\varphi^\infty \big[ e^{-\lambda
\tau} M(\varPhi_\tau^1, \ldots, \varPhi_\tau^N) \big] - M(\varphi)
\end{equation}
for $\varphi \in [0,\infty)^N$ where the infimum is taken over all
(bounded) stopping times $\tau$ of $\varPhi = (\varPhi^1, \ldots,
\varPhi^N)$ and the function $M$ is given by \eqref{4.6} using
\eqref{4.5} above.}

\smallskip

\textbf{Proof.} By It\^o's formula using \eqref{3.3} we get
\begin{equation} \h{3pc} \label{4.8}
e^{-\lambda t} M(\varPhi_t) = M(\varphi) + \int_0^t e^{-\lambda s}
\big( \LL_\varPhi M \m \lambda M)(\varPhi_s)\, ds + N_t
\end{equation}
for $\varphi \in [0,\infty)^N$ \v{-1pt} where $N_t = \sum_{i \in
C_k^n} \sum_{j=1}^k \int_0^t e^{-\lambda s}
M_{\varphi_i}(\varPhi_s)\:\! \mu \:\! \varPhi_s^i\, dB_s^{n_j}$ is a
continuous local martingale for $t \ge 0$. Making use of a
localisation sequence of stopping times for this local martingale if
needed, applying the optional sampling theorem and recalling that
$M$ solves \eqref{4.2}, we find by taking $\EE_\varphi^\infty$ on
both sides in \eqref{4.8} that
\begin{equation} \h{3pc} \label{4.9}
\EE_\varphi^\infty \big[ e^{-\lambda \tau} M(\varPhi_\tau^1,
\ldots, \varPhi_\tau^N) \big] = M(\varphi) + \EE_\varphi^\infty
\Big[ \int_0^\tau\!\! e^{-\lambda t} L(\varPhi_t)\, dt\:\! \Big]
\end{equation}
for all $\varphi \in [0,\infty)^N$ and all (bounded) stopping times
$\tau$ of $\varPhi$. From \eqref{3.4} and \eqref{4.9} using
\eqref{4.3} we see that \eqref{4.7} holds as claimed and the proof
is complete. \hfill $\square$

\smallskip

4.\ From Proposition 2 we see that the optimal stopping problem
\eqref{3.4} is equivalent to the optimal stopping problem defined by
\begin{equation} \h{5pc} \label{4.10}
\check V(\varphi) = \inf_\tau\:\! \EE_\varphi^\infty \big[ e^{-\lambda
\tau} M(\varPhi_\tau^1, \ldots, \varPhi_\tau^N) \big]
\end{equation}
for $\varphi \in [0,\infty)^N$ where the infimum is taken over all
(bounded) stopping times $\tau$ of $\varPhi = (\varPhi^1, \ldots,
\varPhi^N)$ and the function $M$ is given by \eqref{4.6} using
\eqref{4.5} above. The optimal stopping problem \eqref{4.10} is
Mayer formulated. From \eqref{4.7} and \eqref{4.10} we see that
\begin{equation} \h{7.5pc} \label{4.11}
\hat V(\varphi) = \check V(\varphi) - M(\varphi)
\end{equation}
for $\varphi \in [0,\infty)^N$. The Mayer reformulation \eqref{4.10}
has certain advantages that will be exploited in the subsequent
analysis of the optimal stopping problem \eqref{3.4} below.

\section{Hypoellipticity}

In this section we show that the differential operator $\LL_\varPhi$
from \eqref{4.1} satisfies the \emph{H\"ormander condition} (cf.\
\cite{Ho}) and therefore is \emph{hypoelliptic}. This will provide a
needed regularity of the value function $\hat V$ from \eqref{3.4} in
the continuation set of the optimal stopping problem. We also show
that the previous conclusions on the space operator $\LL_\varPhi$
extend to the backward time-space operator $-\partial_t \p
\LL_\varPhi$ which in turn imply that $\varPhi$ is a \emph{strong
Feller} process. This will provide a needed regularity of $\hat V$
at the optimal stopping boundary (between the continuation set and
the stopping set). The regularity of $\hat V$ in the continuation
set and at the optimal stopping boundary will be discussed in
Sections 6 and 7 below.

\smallskip

1.\ \emph{H\"ormander's condition}. To recall what it means that
$\LL_\varPhi$ satisfies the H\"ormander condition, we will connect
to Part 6 of Section 4 in \cite{Pe-4} from where one reads that
$\LL_\varPhi$ can be rewritten as the `sum of squares' as follows
\begin{equation} \h{3pc} \label{5.1}
\LL_\varPhi = \sum_{i=1}^N \mu_i\:\! \partial_{\varphi_i} + \frac{1}
{2} \sum_{i,j=1}^N (\sigma \sigma^t)_{ij}\:\! \partial_{\varphi_i
\varphi_j} = D_0 + \sum_{i=1}^N D_i^2
\end{equation}
where $D_i$ is a first-order differential operator given by
\begin{equation} \h{8pc} \label{5.2}
D_i = \sum_{j=1}^N \beta_{ij}\:\! \partial_{\varphi_j}
\end{equation}
for $0 \le i \le N$ with the coefficients $\beta_{ij}$ expressed
explicitly as \v{-6pt}
\begin{align} \h{5pc} \label{5.3}
&\beta_{0j} = \mu_j - \frac{1}{2} \sum_{k,l=1}^N \sigma_{lk}\;\!
\partial_{\varphi_l} \sigma_{jk}\;\; (1\! \le\! j\! \le\! N) \\[2pt]
\label{5.4} &\beta_{ij} = \tfrac{1}{\sqrt{2}}\, \sigma_{ji}\;\;
(1\! \le\! i \! \le\! N) \; (1\! \le\! j\! \le\! N)\, .
\end{align}
Comparing \eqref{4.1} with \eqref{5.1} we see that the drift vector
$\mu = (\mu_1, \ldots, \mu_N)$ and the diffusion matrix $\sigma =
(\sigma_1; \ldots; \sigma_N)$ of $\varPhi$ (the latter written as a
sequence of its rows) are given by
\begin{align} \h{8pc} \label{5.5}
&\mu_i(\varphi) = \lambda\:\! (1 \p \varphi_i) \\
\label{5.6} &\sigma_i(\varphi) = \mu\:\! \varphi_i\:\! I_i
\end{align}
for $\varphi = (\varphi_1, \ldots, \varphi_N) \in [0,\infty)^N$ and
$1 \le i \le N$ where $I_i = (I_{i1}, \ldots, I_{in},0, \ldots, 0)
\in \{0,1\}^N$ with $I_{ip} = 1_i(p)$ which by definition equals $1$
if $p$ belongs to $\{n_1, \ldots,n_k\}$ for $1 \le i \le N$
identified with $(n_1, \ldots, n_k)$ in $C_k^n$ and $0$ otherwise
(note that we have extended the range of $I_i$ appearing in the text
following \eqref{4.1} above by adding $N \m n$ zeros to make
$\sigma$ a square $N\! \times\! N$ matrix). Identifying $D_i$ with
$(\beta_{i1}, \ldots ,\beta_{iN})$ we see that each $D_i$ may be
viewed as a function from $[0,\infty)^N$ (or its subset) to $\R^N$
defined by $D_i(\varphi) = (\beta_{i1}(\varphi), \ldots
,\beta_{iN}(\varphi))$ for $\varphi \in [0,\infty)^N$ and $0 \le i
\le N$. The \emph{Lie bracket} of $D_i$ and $D_j$ understood as
differential operators is defined by
\begin{equation} \h{6.5pc} \label{5.7}
[D_i,D_j] = D_i D_j \m D_j D_i
\end{equation}
for $0\! \le\! i,j\! \le\! N$. The smallest vector space in $\R^N$
that (i) contains all $D_0,D_1, \ldots, D_N$ understood as vectors
in $\R^N$ and (ii) is closed under the Lie bracket operation
\eqref{5.7} is referred to as the \emph{Lie algebra} generated by
$D_0,D_1, \ldots, D_N$ and is denoted by $Lie\;\! (D_0,D_1, \ldots,
D_N)$. In other words $Lie\:\! (D_0,D_1, \ldots, D_N) =
\text{span}\;\! \{ D_i, [D_i,D_j],[[D_i,D_j],D_k], \ldots\, \vert\,
0 \le i,j,k, \ldots \le N \}$. Note that $Lie\:\! (D_0,D_1, \ldots,
D_N)$ may be viewed as a function from $[0,\infty)^N$ into the
family of linear subspaces of $\R^N$ whose (algebraic) dimensions
could also be strictly smaller than $N$. The H\"ormander condition
states that
\begin{equation} \h{6.5pc} \label{5.8}
\text{dim}\;\! Lie\:\! (D_0,D_1, \ldots, D_N) = N
\end{equation}
on $[0,\infty)^N$ (or its subset). If \eqref{5.8} is satisfied then
$\LL_\varPhi$ is hypoelliptic on $[0,\infty)^N$ (or its subset) by
the H\"ormander theorem (cf.\ \cite[Theorem 1.1]{Ho}).

\medskip

\textbf{Proposition 3.} \emph{The differential operator
$\LL_\varPhi$ from \eqref{4.1} satisfies the H\"ormander condition
\eqref{5.8}. Consequently $\LL_\varPhi$ is hypoelliptic.}

\medskip

\textbf{Proof.} We need to show that \eqref{5.8} holds. For this, we
find by direct calculation using \eqref{5.3} with
\eqref{5.5}+\eqref{5.6} that
\begin{equation} \h{6.5pc} \label{5.9}
\beta_{0j} = a \varphi_j + \lambda\;\; (1\! \le\! j\! \le\! N)
\end{equation}
where we set $a := \lambda \m k \mu^2\!/2$. By \eqref{5.4} and
\eqref{5.6} we see that
\begin{align} \h{6.5pc} \label{5.10}
\beta_{ij} &= \frac{\mu}{\sqrt{2}}\;\! \varphi_j\:\! 1_j(i)\;\;
\text{if}\;\; 1 \le i \le n \\ \notag &= \quad\quad 0\quad\quad
\;\;\; \text{if} \;\;\;\! n < i \le N
\end{align}
for $1 \le j \le N$ where $1_j(i)$ by definition equals $1$ if $i$
belongs to $\{m_1, \ldots, m_k\}$ for $1 \le j \le N$ identified
with $(m_1, \ldots, m_k)$ in $C_k^n$ and $0$ otherwise. Inserting
the right-hand sides of \eqref{5.9} and \eqref{5.10} into
\eqref{5.2} with $i$ \& $j$ swapped we find that
\begin{align} \h{4pc} \label{5.11}
&D_0 = \sum_{i=1}^N (a\:\! \varphi_i \p \lambda)\:\! \partial_
{\varphi_i} \sim \sum_{i=1}^N (\varphi_i \p b)\:\! \partial_
{\varphi_i}\;\; \text{if}\;\; a \ne 0 \\ \notag &\h{8.9pc}\sim
\sum_{i=1}^N \partial_{\varphi_i}\;\; \text{if}\;\; a=0 \\
\label{5.12} &D_j = \sum_{i=1}^N \frac{\mu}{\sqrt{2}}\,
\varphi_i\:\! I_{ij}\:\! \partial_{\varphi_i} \sim \sum_{i=1}^N
\varphi_i\:\! I_{ij}\:\! \partial_{\varphi_i}\;\; (1\! \le\! j
\! \le\! N)
\end{align}
where in \eqref{5.11} we set $b := \lambda/a$ and $L \sim R$ by
definition means that $R$ is a constant multiple of $L$ (\:\!making
$R$ equivalent to $L$ when searching for the Lie algebra generated
by a set containing $L$) and $I_{ij}$ in \eqref{5.12} equals $1$ if
$j$ belongs to $\{n_1, \ldots, n_k\}$ for $1 \le i \le N$ identified
with $(n_1, \ldots\, n_k)$ in $C_k^n$ and $0$ otherwise.

Calculating (iterated) Lie brackets of $D_0,D_1, \ldots, D_N$ from
\eqref{5.11} and \eqref{5.12} one arrives at the following recipe
for verifying the H\"ormander condition \eqref{5.8}. For this, let
any $1 \le j \le N$ with $j = (n_1, \ldots,n_k) \in C_k^n$ be given
and fixed. Using \eqref{5.11} and \eqref{5.12} we find by direct
calculation that
\begin{align} \h{4pc} \label{5.13}
[D_0,D_{n_1}] &\sim \sum_{i=1}^N b\:\! I_{i\:\!n_1} \partial_{\varphi_i}
\sim \sum_{i=1}^N I_{i\:\!n_1} \partial_{\varphi_i}\;\; \text{if}\;\;
a \ne 0 \\ \notag &= \sum_{i=1}^N I_{i\:\!n_1} \partial_{\varphi_i}
\;\; \text{if}\;\; a=0
\end{align}
for any $1 \le n_1 \le n$ and thus the specified one in forming $j$
as well. Using \eqref{5.12} and \eqref{5.13}we find by direct
calculation that
\begin{equation} \h{6pc} \label{5.14}
[\:\![D_0,D_{n_1}],D_{n_2}] \sim \sum_{i=1}^N I_{i\:\!n_1} I_{i\:\!n_2}
\:\! \partial_{\varphi_i}
\end{equation}
for any $1 \le n_1 < n_2 \le n$ and thus the specified ones in
forming $j$ as well. Continuing this calculation by induction we
find that
\begin{equation} \h{4pc} \label{5.15}
[\:\![\:\![D_0,D_{n_1}],D_{n_2}], \ldots, D_{n_k}] \sim \sum_{i=1}^N
I_{i\:\!n_1} I_{i\:\!n_2} \ldots I_{i\:\!n_k} \:\! \partial_{\varphi_i}
\end{equation}
for any $1 \le n_1 < n_2 < \ldots < n_k \le n$ and thus the
specified ones in forming $j$ as well. By definition of $I$ recalled
following \eqref{5.12} above we see that
\begin{align} \h{6pc} \label{5.16}
I_{i\:\!n_1} I_{i\:\!n_2} \ldots I_{i\:\!n_k} &= 1\;\; \text{if}
\;\; i = j \\[3pt] \notag &= 0\;\; \text{if}\;\; i \ne j\, .
\end{align}
Combining \eqref{5.15} and \eqref{5.16} we find that
\begin{equation} \h{6pc} \label{5.17}
[\:\![\:\![D_0,D_{n_1}],D_{n_2}], \ldots, D_{n_k}] \sim
\partial_{\varphi_j}\, .
\end{equation}
Since $1 \le j \le N$ was arbitrary this shows that the iterated Lie
brackets \eqref{5.17} span the entire $\R^N$ so that the H\"ormander
condition \eqref{5.8} is satisfied as claimed. The final claim
follows by the H\"ormander theorem as recalled following \eqref{5.8}
above. This completes the proof. \hfill $\square$

\medskip

2.\ \emph{Parabolic H\"ormander's condition}. To recall what it
means that $-\partial_t \p \LL_\varPhi$ satisfies the (parabolic)
H\"ormander condition, we will connect to Part 3 of Section 5 in
\cite{Pe-4} from where one reads that $-\partial_t \p \LL_\varPhi$
can be rewritten as the `sum of squares' as follows
\begin{equation} \h{2pc} \label{5.18}
-\partial_t \p \LL_\varPhi = -\partial_{\varphi_0} + \sum_{i=1}^N
\mu_i\:\! \partial_{\varphi_i} + \frac{1}{2} \sum_{i,j=1}^N (\sigma
\sigma^t)_{ij}\:\! \partial_{\varphi_i \varphi_j} = \bar D_0 +
\sum_{i=1}^N \bar D_i^2
\end{equation}
where $\bar D_i$ is a first-order differential operator given by
\begin{equation} \h{8pc} \label{5.19}
\bar D_i = \sum_{j=0}^N \beta_{ij}\:\! \partial_{\varphi_j}
\end{equation}
for $0 \le i \le N$ with the coefficients $\beta_{ij}$ expressed
explicitly as
\begin{equation} \h{6pc} \label{5.20}
\beta_{00} = -1\;\; \&\;\; \beta_{i0} = 0\;\; (1\! \le\! i\! \le\! N)
\end{equation}
in addition to \eqref{5.3} and \eqref{5.4} above. Viewing $\bar D_i$
as functions from $[0,\infty)^{N+1}$ (or its subset) to $\R^{N+1}$
this amounts to setting
\begin{equation} \h{3pc} \label{5.21}
\bar D_0 = (-1,\beta_{01}, \ldots\, \beta_{0N})\;\; \&\;\; \bar D_i
= (0,\beta_{i1}, \ldots, \beta_{iN})\;\; (1\! \le\! i\! \le\! N)\, .
\end{equation}
Note that $Lie\:\! (\bar D_0,\bar D_1, \ldots, \bar D_N)$ can be
viewed as a function from $[0,\infty)^{N+1}$ into the family of
linear subspaces of $\R^{N+1}$ whose (algebraic) dimensions could
also be strictly smaller than $N \p 1$. The (parabolic) H\"ormander
condition states that
\begin{equation} \h{6pc} \label{5.22}
\text{dim}\;\! Lie\:\! (\bar D_0,\bar D_1, \ldots, \bar D_N) = N \p 1
\end{equation}
on $[0,\infty)^{N+1}$ (or its subset). If \eqref{5.22} is satisfied
then $-\partial_t \p \LL_\varPhi$ is hypoelliptic on
$[0,\infty)^{N+1}$ (or its subset) by the H\"ormander theorem (cf.\
\cite[Theorem 1.1]{Ho}).

\medskip

\textbf{Proposition 4.} \emph{The backward time-space differential
operator $-\partial_t + \LL_\varPhi$ of $\varPhi$ satisfies the
parabolic H\"ormander condition \eqref{5.22}. Consequently $\varPhi$
is a strong Feller process.}

\medskip

\textbf{Proof.} 1.\ We need to show that \eqref{5.22} holds. For
this, note from \eqref{5.12} that with a slight abuse of notation we
have
\begin{equation} \h{5pc} \label{5.23}
\bar D_0 = (-1,D_0)\;\; \&\;\; \bar D_i = (0,D_i)\;\; (1\! \le\! i
\! \le\! N)
\end{equation}
where $D_0$ and $D_i$ are defined in \eqref{5.2} with
\eqref{5.3}+\eqref{5.4} above followed by the more explicit
expressions in \eqref{5.11}+\eqref{5.12} above. Since the first
coordinates $-1$ and $0$ in \eqref{5.23} are constants, and
$\partial_{\varphi_0} = \partial_t$ plays no role because $\mu$ and
$\sigma$ are time independent, we see that calculating the
(iterated) Lie brackets of $\bar D_0, \bar D_1, \ldots, \bar D_N$
reduces to calculating the (iterated) Lie brackets of $D_0, D_1,
\ldots, D_N$. More specifically, letting any $1 \le i \le N$ with $i
= (n_1, \ldots,n_k) \in C_k^n$ be given and fixed, and proceeding in
the same way as in the proof of Proposition 3, we find by
\eqref{5.17} above that
\begin{equation} \h{3pc} \label{5.24}
[\:\![\:\![\bar D_0,\bar D_{n_1}],\bar D_{n_2}], \ldots, \bar
D_{n_k}] \sim [\:\![\:\! [D_0,D_{n_1}],D_{n_2}], \ldots, D_{n_k}]
\sim \partial_{\varphi_i}\, .
\end{equation}
Viewed as a function from $[0,\infty)^{N+1}$ to $\R^{N+1}$ we see
that the differential operator \eqref{5.24} is identified with $e_i
= (0, \ldots,0,1,0, \ldots, 0) \in \R^{N+1}$ where the number $1$ is
placed on the $i\p 1$ coordinate. Enlarging \eqref{5.24} by
\begin{equation} \h{6pc} \label{5.25}
\bar D_0 = (-1,D_0) = -\partial_{\varphi_0} + \sum_{i=1}^N (a \varphi_i
\p \lambda)\:\! \partial_ {\varphi_i}
\end{equation}
which viewed as a function from $[0,\infty)^{N+1}$ to $\R^{N+1}$ is
identified with $e_0 = (-1, a\varphi_1 \p \lambda, \ldots, a
\varphi_N \p \lambda)$ for $(\varphi_1, \ldots, \varphi_N) \in
[0,\infty)^N$, we see that the vectors $e_i$ for $0 \le i \le N$ are
linearly independent in $\R^{N+1}$ and clearly span the entire
$\R^{N+1}$ so that the parabolic H\"ormander condition \eqref{5.22}
is satisfied as claimed.

\smallskip

2.\ To see that $\varPhi$ is a strong Feller process, let any
bounded measurable function $f : [0,\infty)^N \rightarrow \R$ be
given and fixed. Then by Corollary 9 in \cite{Pe-4} we know that
\begin{equation} \h{7pc} \label{5.26}
\partial_t \PP_{\!t} f \buildrel w \over = \LL_\varPhi \PP_{\!t} f
\;\;\; \text{in}\;\; (0,\infty)\! \times\! [0,\infty)^N
\end{equation}
where $\PP_{\!t} f(\varphi) = \EE_\varphi^\infty f(\varPhi_t)$ for
$t \ge 0$ and $\varphi \in [0,\infty)^N$. Since $-\partial_t +
\LL_\varPhi$ satisfies the H\"orman- der condition \eqref{5.22} we
know that $-\partial_t + \LL_\varPhi$ is hypoelliptic by the
H\"ormander theorem as recalled following \eqref{5.22} above. It
follows therefore that the weak solution $(t,\varphi) \mapsto
\PP_{\!t} f(\varphi)$ in \eqref{5.26} is $C^\infty$ on $(0,\infty)\!
\times\! [0,\infty)^N$. In particular, the mapping $\varphi \mapsto
\PP_{\!t} f(\varphi)$ is $C^\infty$ and therefore continuous as well
on $[0,\infty)^N$ for every $t>0$ given and fixed. This shows that
$\varPhi$ is a strong Feller process as claimed and the proof is
complete. \hfill $\square$

\section{Properties of the optimal stopping boundary}

In this section we establish the existence of an optimal stopping
time in the problem \eqref{3.4} and derive basic properties of the
optimal stopping boundary.

1.\ Looking at \eqref{3.4} we may conclude that the (candidate)
continuation and stopping sets in this problem need to be defined as
follows
\begin{align} \h{7pc} \label{6.1}
C = \{\, \varphi \in [0,\infty)^N\; \vert\; \hat V(\varphi)<0\, \}
\\[2pt] \label{6.2} D = \{\, \varphi \in [0,\infty)^N\; \vert\;
\hat V(\varphi) = 0\, \}
\end{align}
respectively. Since the integral in the optimal stopping problem
\eqref{3.4} is uniformly bounded from below by $-1/c$, it follows by
\cite[Theorem 6(2), p.\ 137]{Sh-2} that the first entry time of the
process $\varPhi$ into the closed set $D$ defined by
\begin{equation} \h{7pc} \label{6.3}
\tau_D = \inf\, \{\, t \ge 0\; \vert\; \varPhi_t \in D\, \}
\end{equation}
is optimal in \eqref{4.3}. Within this general existence result one
formally allows that the optimal stopping time $\tau_D$ takes the
value $\infty$ as well, however, we will now show that this is not
the case in the present problem.

\medskip

\textbf{Proposition 5.} \emph{We have
\begin{equation} \h{8.5pc} \label{6.4}
\PP_{\!\varphi}^\infty(\tau_D < \infty) = 1
\end{equation}
for all $\varphi \in [0,\infty)^N$.}

\medskip

\textbf{Proof.} Since $\tau_D$ is optimal in \eqref{3.4} it follows
from the result and proof of Proposition 1 above that $\tau_D$ is
optimal in \eqref{2.8} as well. The result of Lemma 1 in \cite{JP}
identifies the Radon-Nikodym derivative corresponding to the measure
change from $\PP_{\!\pi}$ to $\PP^\infty$ to be
\begin{equation} \h{8.5pc} \label{6.5}
\frac{d\PP_{\!\pi,\tau}}{d\PP_{\!\tau}^\infty} = e^{-\lambda \tau}
\:\! \frac {1 \m \pi}{1 \m \Pi_\tau}
\end{equation}
for all stopping times $\tau$ of $X$ and all $\pi \in [0,1)$, where
$\PP_{\!\tau}^\infty$ and $\PP_{\!\pi,\tau}$ denote the restrictions
of measures $\PP^\infty$ and $\PP_{\!\pi}$ to $\cF_\tau^X$ for $\pi
\in [0,1)$ respectively. Using \eqref{6.5} we recognise $e^{\lambda
t} (1 - \Pi_t)$ as a constant multiple of the Radon-Nikodym
derivative $d\PP_{\!t}^\infty/d\PP_{\!\pi,t}$ and hence the process
is a martingale under $\PP_{\!\pi}$ for $t \ge 0$ whenever $\pi \in
[0,1)$ is given and fixed. Moreover, from the fact that the
probability measures $\PP_{\!\pi}$ and $\PP^\infty$ are singular, we
can conclude that
\begin{equation} \h{8.5pc} \label{6.6}
e^{\lambda t} (1 \m \Pi_t) \rightarrow 0
\end{equation}
and hence $\Pi_t \rightarrow 1$ both almost surely under
$\PP_{\!\pi}$ as $t \rightarrow \infty$ for $\pi \in [0,1)$ (cf.\
Theorem 2 in \cite[p.\ 527]{Sh-3}). Using the latter fact on the
right-hand side of \eqref{2.8} with $\tau_D$ in place of $\tau$ we
see that $\PP_{\!\pi}(\tau_D < \infty) = 1$ since otherwise $1 \m
\pi \ge V(\pi) = \infty$ for $\pi \in [0,1)$ which is a
contradiction. Since the set $\{ \tau_D < \infty \}$ belongs to
${\cal F}_{\tau_D}^X$ and by \eqref{6.5} above the probability
measures $\PP_{\!\pi}$ and $\PP^\infty$ restricted to ${\cal
F}_{\tau_D}^X$ are equivalent (i.e.\ $\PP_{\!\pi}(F) = 0$ if and
only if $\PP^\infty(F) = 0$ for $F \in {\cal F}_{\tau_D}^X$) for
$\pi \in [0,1)$, it follows that $\PP_{\!\varphi}^\infty ( \tau_D <
\infty ) = 1$ for all $\varphi \in [0,\infty)^N$ as claimed. \hfill
$\square$

\medskip

2.\ The topological boundary between the sets $C$ and $D$ in
$[0,\infty)^N$ is referred to as the \emph{optimal stopping
boundary} in the problem \eqref{3.4}. We will denote it by $\partial
C$ although we could also use $\partial D$ without altering its
meaning. To derive an upper bound on the size of $\partial C$,
recall that the optimal stopping boundary/point $\varphi_*$ in the
one-dimensional problem \eqref{3.4} (i.e.\ when $N=1$) can be
characterised as a unique solution to
\begin{equation} \h{3pc} \label{6.7}
\frac{e^{\kappa(1+\varphi_*)/\varphi_*}}{\varphi_*^\kappa} \int_0^{
\varphi_*/(1 + \varphi_*)}\!\! \frac{u^{\kappa-1}} {(1 \m u)^{\kappa
+2}} \, e^{-\kappa/u}\, du = \frac{\mu^2}{2c}
\end{equation}
on $(\lambda/c,\infty)$ where we set $\kappa = 2 \lambda/\mu^2$
(cf.\ \cite[Section 5]{EP} \& \cite[Section 22]{PS}). Note from
\eqref{6.7} that we have $\varphi_* =
\varphi_*(\lambda/\mu^2,\lambda/c)$ and set
\begin{equation} \h{8pc} \label{6.8}
\varphi_i^* := \varphi_* \big( \tfrac{\lambda}{k \mu^2},\tfrac{\lambda}
{p_i c} \big)
\end{equation}
for $\lambda>0$, $\mu \in \R$,$c>0$ and $p_i \in [0,1]$ for $1 \le i
\le N$ with $\sum_{i=1}^N p_i = 1$. Recall that $\varphi_i^* \in
(\lambda/(p_i c),\infty)$ for $1 \le i \le N$. We can now expose
basic properties of the the value function and the
continuation/stopping set in the problem \eqref{3.4} as follows.

\medskip

\textbf{Proposition 6.} \emph{\begin{align} \h{0pc} \label{6.9}
&\h{-5pt}\text{The value function}\;\; \hat V\;\; \text{is concave
and continuous on}\;\; [0,\infty)^N\, . \\[3pt] \label{6.10}
&\h{-5pt}\text{If}\;\; \varphi_1 \le \psi_1, \ldots, \varphi_N \le
\psi_N\;\; \text{then}\;\; \hat V(\varphi_1, \ldots, \varphi_N)
\le \hat V(\psi_1, \ldots \psi_N)\, . \\[3pt] \label{6.11}
&\h{-5pt} \text{If}\;\; (\varphi_1, \ldots, \varphi_N) \in
D\;\; \text{and}\;\; \psi_1 \ge \varphi_1, \ldots, \psi_N \ge
\varphi_N\;\; \text{then}\;\; (\psi_1, \ldots \psi_N) \in D\, .
\\[3pt] \label{6.12} &\h{-5pt}\text{The stopping set}\;\; D\;\;
\text{is convex}\;\; \text{and the polytope}\;\; \{\, (\varphi_1,
\ldots \varphi_N) \in [0,\infty)^N\; \vert\; \\[-3pt] \notag
&\h{-3pt} \textstyle \sum_{i=1}^N \varphi_i/\varphi_i^* \m 1
\ge 0\, \} \;\; \text{is contained in}\;\; D\, . \\[3pt]
\label{6.13} &\h{-5pt} \text{The simplex} \;\; \{\, (\varphi_1,
\ldots, \varphi_N) \in [0,\infty)^N\; \vert\; \textstyle
\sum_{i=1}^N p_i \varphi_i \m \lambda/c < 0\, \}\;\; \text{is
contained} \\[-3pt] \notag &\h{-5pt}\text{in the continuation
set}\;\; C\, .
\end{align}}

\textbf{Proof.} \eqref{6.9}: Combining the fact that the Markovian
flow \eqref{2.12} is linear as a function of its initial point with
the fact that the integral in \eqref{3.4} is a linear function of
its argument, and using that the infimum of a convex combination is
larger than the convex combination of the infima, we find that $\hat
V$ is concave on $[0,\infty)^N$ as claimed. Hence we can also
conclude that $\hat V$ is continuous on the open set $(0,\infty)^N$.
To see that $\hat V$ is continuous at the boundary points of
$[0,\infty)^N$ we may recall the well-known (and easily verified)
fact that the concave function $\hat V$ is lower semicontinuous on
the closed and convex set $[0,\infty)^N$. Moreover, recalling that
\eqref{2.12} defines a Markovian functional of the initial point
$\varPhi_0^i := \varphi_i$ in $[0,\infty)$ of the process
$\varPhi^i$ for $1 \le i \le N$, we see that the expectation in
\eqref{4.10} defines a continuous function of the initial point
$\varphi = (\varphi_1, \ldots, \varphi_N)$ of the process $\varPhi =
(\varPhi^1, \ldots, \varPhi^N)$ for every (bounded) stopping time
$\tau$ of $\varPhi$ given and fixed. Taking the infimum over all
(bounded) stopping times $\tau$ of $\varPhi$ we can thus conclude
from \eqref{4.11} that the value function $\hat V$ is upper
semicontinuous on $[0,\infty)^N$. Being also lower semicontinuous it
follows that $\hat V$ is continuous on $[0,\infty)^N$ as claimed.

\smallskip

\eqref{6.10}: This is a direct consequence of the fact that the
Markovian flow \eqref{2.12} is increasing as a function of its
initial point being used in \eqref{3.4} above.


\begin{figure}[!t] 
\begin{center}
\includegraphics[scale=0.9]{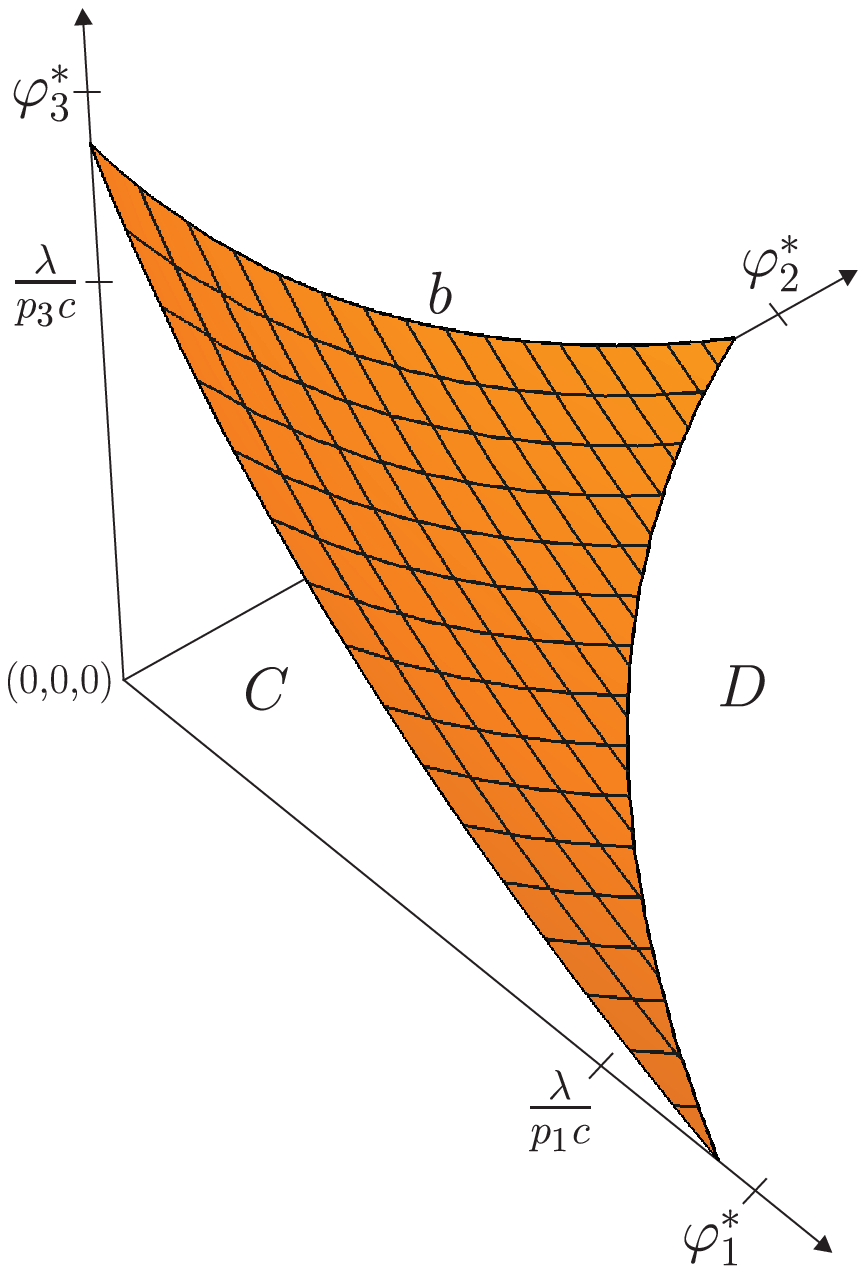}
\end{center}

{\par \leftskip=2cm \rightskip=2cm \small \noindent

\textbf{Figure 1.} The optimal stopping surface $b$ in the problem
\eqref{3.4} when $n=3$ and $k=2$ with $\mu=\lambda=c=1$ and $p_{1,2}
= p_{1,3} = p_{2,3} = 1/3$.

\par} 

\end{figure}



\smallskip

\eqref{6.11}: By \eqref{6.10} we have $\hat V(\varphi_1, \ldots,
\varphi_N) \le \hat V(\psi_1, \ldots, \psi_N) \le 0$ so that
$(\varphi_1, \ldots, \varphi_N) \in D$ i.e.\ $\hat V(\varphi_1,
\ldots, \varphi_N) = 0$ implies that $\hat V(\psi_1, \ldots, \psi_N)
= 0$ i.e.\ $(\psi_1, \ldots, \psi_N) \in D$ as claimed.

\eqref{6.12}: To see that $D$ is convex, take any $\varphi$ and
$\psi$ from $D$ and note by \eqref{6.9} that $0 \ge$ $\hat V(\alpha
\varphi \p (1 \m \alpha) \psi) \ge \alpha \hat V(\varphi) \p (1 \m
\alpha) \hat V(\psi) = 0$ so that $\hat V(\alpha \varphi \p (1 \m
\alpha) \psi) = 0$ i.e.\ $\alpha \varphi \p (1 \m \alpha) \psi$ $\in
D$ for every $\alpha \in [0,1]$ as claimed. To see that the polytope
is contained in $D$, note that pulling $p_j$ in front of the infimum
in \eqref{3.4} with any $1 \le j \le N$ given and fixed shows that
the point $(0, \dots,0,\varphi_j^*,0, \ldots,0)$ belongs to $D$
because $\varphi_j^*$ as defined in \eqref{6.8} above (\:\!with $j$
in place of $i$) is \v{-1pt} an optimal stopping point in the
one-dimensional problem obtained by removing the non-negative term
$\sum_{i=1,i \ne j}^N (p_i/p_j)\;\! \varPhi_t^i$ from the integral
with respect to time in \eqref{3.4} with $p_j$ in front of \v{-1pt}
the infimum (\:\!note that the appearance of $k$ in \eqref{6.8}
follows from the fact that $(1/\sqrt{k}) \sum_{j=1}^k B_t^{n_j}$ is
a standard Brownian motion for $t \ge 0$). It follows therefore by
\eqref{6.11} that the set $\{\;\! (0, \dots,0,\varphi_j,0,
\ldots,0)\; \vert\; \varphi_j \ge \varphi_j^*\, \}$ is contained in
$D$ for every $1 \le j \le N$. But then the entire polytope is
contained in $D$ due to its convexity.

\smallskip

\eqref{6.13}: Taking any point $\varphi$ from the (open) simplex and
replacing $\tau$ in \eqref{3.4} by the first exit time of $\varPhi$
from a sufficiently small ball around $\varphi$ that is strictly
contained in the simplex, we see that the integrand in \eqref{3.4}
remains strictly negative so that $\hat V$ takes a strictly negative
value at $\varphi$ itself, showing that $\varphi$ belongs to the
continuation set $C$ as claimed. \hfill $\square$

\medskip

3.\ From the results of Proposition 6 we see that the stopping set
in the problem \eqref{3.4} can be described as follows
\begin{equation} \h{4pc} \label{6.14}
D = \{\, (\varphi_1, \ldots, \varphi_N) \in [0,\infty)^N\; \vert\;
\varphi_N \ge b(\varphi_1, \ldots, \varphi_{N-1})\, \}
\end{equation}
where $b : [0,\infty)^{N-1} \rightarrow [0,\infty)$ is a convex,
continuous, decreasing function (\:\!in the sense that $b(\varphi_1,
\ldots, \varphi_{N-1}) \ge b(\psi_1, \ldots, \psi_{N-1})$ whenever
$\varphi_1 \le \psi_1, \ldots, \varphi_{N-1} \le \psi_{N-1}$)
satisfying
\begin{equation} \h{3pc} \label{6.15}
-\sum_{i=1}^{N-1} \frac{p_i}{p_N}\, \varphi_i + \frac{\lambda}{p_N c}
\le b(\varphi_1, \ldots, \varphi_{N-1}) \le - \sum_{i=1}^{N-1}
\frac{\varphi_N^*}{\varphi_i^*}\, \varphi_i + \varphi_N^*
\end{equation}
for $\varphi_1 \in [0,\lambda/(p_1 c)], \ldots, \varphi_{N-1} \in
[0,\lambda/(p_{N-1} c)]$ and $\varphi_1 \in [0,\varphi_1^*], \ldots,
\varphi_{N-1} \in [0,\varphi_{N-1}^*]$ in the first and second
inequality respectively (see Figure 1). Note that the optimal
stopping boundary in the problem \eqref{3.4} can be described as
follows
\begin{equation} \h{4pc} \label{6.16}
\partial C = \{\, (\varphi_1, \ldots, \varphi_N) \in [0,\infty)^N\;
\vert\; \varphi_N = b(\varphi_1, \ldots, \varphi_{N-1})\, \}\, .
\end{equation}
We address the question of characterising/determining $b$ in the
remaining two sections. To this end we conclude this section by
establishing a key regularity result of $\partial C$ for $D$.

\medskip

4.\ Recall that a point $\varphi \in [0,\infty)^N$ is said to be
\emph{probabilistically regular} for $D$ if
\begin{equation} \h{10pc} \label{6.17}
\PP_{\!\varphi}^\infty (\sigma_D\! =\! 0) = 1
\end{equation}
where $\sigma_D$ is the first hitting time of $\varPhi$ to $D$
defined by $\sigma_D = \inf\, \{\, t > 0\; \vert\; \varPhi_t \in D\,
\}$ (see Sections 2 and 3 in \cite{DePe} for fuller details). If
every point at $\partial C$ is probabilistically regular for $D$ we
say that $\partial C$ is probabilistically regular for $D$. We now
show that this is the case in the optimal stopping problem
\eqref{3.4}.

\medskip

\textbf{Proposition 7.} \emph{The optimal stopping boundary
$\partial C$ is probabilistically regular for the stopping set $D$
in the problem \eqref{3.4}.}

\medskip

\textbf{Proof.} Let any point $\varphi = (\varphi_1, \ldots,
\varphi_N)$ at $\partial C$ be given and fixed. We need to show that
\eqref{6.17} holds. For this, note that
\begin{align} \h{1pc} \label{6.18}
\PP_{\!\varphi}^\infty \big(\sigma_D\! =\! 0) &= \PP_{\!\varphi}^
\infty (\cap_{n=1}^\infty \cup_{t \in (0,1/n)} \{ \varPhi_t\! \in
\! D \} \big) = \lim_{n \rightarrow \infty} \PP_{\!\varphi}^\infty
\big(\cup_{t \in (0,1/n)} \{ \varPhi_t\! \in\! D \} \big) \\[0pt]
\notag &\ge \lim_{n \rightarrow \infty} \sup_{t \in (0,1/n)} \PP_
{\!\varphi}^ \infty (\varPhi_t\! \in\! D) = \limsup_{t \downarrow
0}\;\! \PP_{\! \varphi}^ \infty (\varPhi_t\! \in\! D) \\[2pt] \notag
&\ge \limsup_{t \downarrow 0}\;\! \PP_{\!\varphi}^\infty \big(
\varPhi_t^1 \ge \varphi_1, \ldots, \varPhi_t^N \ge \varphi_N \big)
\\[1.5pt] \notag &\ge \limsup_{t \downarrow 0}\;\! \PP^\infty
\big( L_t^1 \ge 1, \ldots, L_t^N \ge 1 \big) \\[-2pt]
\notag &= \limsup_{t \downarrow 0}\;\! \PP^\infty \Big( \textstyle
\sum_{j=1}^k X_t^{n_j^1} \ge k\;\! \tfrac {\mu}{2}\:\! t, \ldots,
\sum_{j=1}^k X_t^{n_j^N}\! \ge k\;\! \tfrac{\mu}{2}\:\! t \Big)
\\[1.5pt] \notag &\ge \limsup_{t \downarrow 0}\;\! \PP^\infty
\big( X_t^1 \ge \tfrac{\mu}{2}\:\! t, \ldots, X_t^n \ge \tfrac
{\mu}{2} \:\! t \big) \\[1pt] \notag &= \limsup_{t \downarrow
0}\;\! \PP \big( B_t^1 \ge \tfrac{\mu}{2}\:\! t, \ldots, B_t^n
\ge \tfrac{\mu}{2} \:\! t \big) \\ \notag &= \limsup_{t \downarrow
0}\;\! \big( \PP \big( B_t^1 \ge \tfrac{\mu}{2}\:\! t \big) \big)^n
= \lim_{t \downarrow 0}\;\! \big( \PP \big( B_1^1 \ge \tfrac{\mu}
{2}\;\! \sqrt{t} \big) \big)^n \\[1pt] \notag &= \big( \PP (B_1^1\!
>\! 0) \big)^n = 1/2^n > 0
\end{align}
where in the second inequality we use that $[\varphi_1,\infty)\!
\times \ldots \times [\varphi_N,\infty) \subseteq D$ by \eqref{6.11}
above, in the third inequality we use \eqref{2.12} above, and in the
fourth equality we use \eqref{2.11} above (upon assuming that $\mu >
0$ without loss of generality). From \eqref{6.18} we see that
$\PP_{\!\varphi}^\infty \big(\sigma_D\! =\! 0) > 0$. As the latter
probability can only be either zero or one by the Blumenthal 0-1 law
(cf.\ \cite[p.\ 30]{BG}), it follows that \eqref{6.17} holds as
claimed. \hfill $\square$

\medskip

\textbf{Remark 8.} Note that if we replace all inequalities under
the probability measures in \eqref{6.18} by strict inequalities,
then the same proof shows that $\partial C$ is probabilistically
regular for the interior of the stopping set $D$ in the problem
\eqref{3.4}. Although this stronger probabilistic regularity plays
no role in the present setting because the process $\varPhi$ is
strong Feller by Proposition 4 above, this observation may be useful
in the settings where the process $\varPhi$ is only known to be
strong Markov (see \cite{DePe} for fuller details).

\section{Free-boundary problem}

In this section we derive a free-boundary problem that stands in
one-to-one correspondence with the optimal stopping problem
\eqref{3.4}. Using the results derived in the previous sections we
show that the value function $\hat V$ from \eqref{3.4} and the
optimal stopping boundary $b$ from \eqref{6.16} solve the
free-boundary problem. This establishes the existence of a solution
to the free-boundary problem. Its uniqueness in a natural class of
functions will follow from a more general uniqueness result that
will be established in Section 8 below. This will also yield an
explicit  integral representation of the value function $\hat V$
expressed in terms of the optimal stopping boundary $b$.

\medskip

1.\ Consider the optimal stopping problem \eqref{3.4} where the
Markov process $\varPhi = (\varPhi^1, \ldots, \varPhi^N)$ solves the
system of stochastic differential equations \eqref{3.3} driven by a
standard $n$\!-dimensional Brownian motion $B = (B^1, \ldots, B^n)$
under the probability measure $\PP^\infty$. Recall that the
infinitesimal generator of $\varPhi$ is the second-order
hypoelliptic differential operator $\LL_\varPhi$ given in
\eqref{4.1} above (cf.\ Proposition 3). Looking at \eqref{3.4} and
relying on other properties of $\hat V$ and $b$ derived above, we
are naturally led to formulate the following free-boundary problem
for finding $\hat V$ and $b$ :
\begin{align} \h{4pc} \label{7.1}
&\LL_\varPhi \hat V \m \lambda \hat V = -L\;\;\; \text{in}\;\; C
\\ \label{7.2} &\hat V(\varphi) = 0\;\;\; \text{for}\;\; \varphi
\in D\;\;\; \text{(instantaneous stopping)} \\ \label{7.3} &\hat
V_{\varphi_i}(\varphi) = 0\;\;\; \text{for}\;\; \varphi \in
\partial C\;\; \text{and} \;\; i=1, \ldots, N\;\;\; \text{(smooth
fit)}
\end{align}
where $L$ is defined in \eqref{4.3} above, $C$ is the (continuation)
set from \eqref{6.1} above, $D$ is the (stopping) set from
\eqref{6.2}+\eqref{6.14} above, and $\partial C$ is the (optimal
stopping) boundary between the sets $C$ and $D$ from \eqref{6.16}
above.

\medskip

2.\ To formulate the existence and uniqueness result for the
free-boundary problem \eqref{7.1}-\eqref{7.3}, we let $\cal C$
denote the class of functions $(U,a)$ such that
\begin{align} \h{0pc} \label{7.4}
&U\;\; \text{belongs to}\;\; C^2(C_a) \cap C^1(\bar C_a)\;\;
\text{and}\;\; \text{is continuous \& bounded on}\;\; [0,\infty)^N
\\[2pt] \label{7.5} &a\;\; \text{is continuous \& decreasing on}\;\;
[0,\infty)^{N-1}\; (\:\!\text{in the sense that}\;\; a(\varphi_1,
\ldots, \varphi_{N-1}) \\[-2pt] \notag &\ge a(\psi_1, \ldots,
\psi_{N-1})\;\; \text{whenever}\;\; \varphi_1 \le \psi_1, \ldots,
\varphi_{N-1} \le \psi_{N-1})\;\; \text{and satisfies}\;\; \\[-3pt]
\notag &\textstyle \sum_{i=1}^{N-1} p_i \varphi_i + p_N a(\varphi_1,
\ldots, \varphi_{N-1}) \m \lambda/c \ge 0\;\; \text{for}\;\;
(\varphi_1, \ldots, \varphi_{N-1}) \in [0,\infty)^{N-1}
\end{align}
where we set $C_a = \{\;\! (\varphi_1, \ldots, \varphi_N) \in
[0,\infty)^N\; \vert\; \varphi_N < a(\varphi_1, \ldots,
\varphi_{N-1})\;\! \}$ and $\bar C_a = \{\, (\varphi_1, \ldots,
\varphi_N) \\ \in [0,\infty)^N\; \vert\; \varphi_N \le a(\varphi_1,
\ldots, \varphi_{N-1})$ with $(\varphi_1, \ldots, \varphi_{N-1})$
belonging to the closure of the set $\{\, (\psi_1, \ldots,
\psi_{N-1}) \; \vert \; a(\psi_1, \ldots, \psi_{N-1})>0\, \}$ in
$[0,\infty)^{N-1}\, \}$.

\medskip

\textbf{Theorem 9.} \emph{The free-boundary problem
\eqref{7.1}-\eqref{7.3} has a unique solution $(\hat V,b)$ in the
class $\cal C$ where $\hat V$ is given in \eqref{3.4} and $b$ is
given in \eqref{6.16} above.}

\medskip

\textbf{Proof.} We first show that the pair $(\hat V,b)$ belongs to
the class $\cal C$ and solves the free-boundary problem
\eqref{7.1}-\eqref{7.3}. For this, note that the optimal stopping
problem \eqref{3.4} is Lagrange formulated and recall that the
infimum in \eqref{3.4} is attained at $\tau_D$ from \eqref{6.3}
above. It follows therefore by the result of Corollary 6 in
\cite{Pe-4} that $\hat V$ from \eqref{3.4} is a weak solution to the
equation \eqref{7.1} in the sense of Schwartz distributions. By the
result of Proposition 3 above we know that the differential operator
$\LL_\varPhi$ is hypoelliptic and hence $\LL_\varPhi \m \lambda I$
is hypoelliptic too. It follows therefore that the weak solution
$\hat V$ to \eqref{7.1} belongs to $C^\infty$ on $C$ (cf.\
\cite[Corollary 8]{Pe-4}). This shows that $\hat V$ belongs to
$C^2(C)$ and satisfies \eqref{7.1} as claimed. Moreover, from
\eqref{6.9} we know that $\hat V$ is continuous on $[0,\infty)^N$
and from \eqref{3.4} we readily find that
\begin{equation} \h{9pc} \label{7.6}
- \frac{1}{c} \le \hat V(\varphi) \le 0
\end{equation}
for all $\varphi \in [0,\infty)^N$. Furthermore, recall from
Proposition 4 above that the process $\varPhi = (\varPhi^1, \ldots,
\varPhi^N)$ is strong Feller while by Proposition 7 we know that
$\partial C$ is probabilistically regular for $D$. Finally, from
\eqref{2.12} we see that the process $\varPhi$ can be realised as a
continuously differentiable stochastic flow of its initial point so
that the integrability conditions of Theorem 8 in \cite{DePe} are
satisfied. Recalling that $\hat V$ satisfies \eqref{7.2}, and
applying the result of that theorem, we can conclude that
\begin{equation} \h{4.5pc} \label{7.7}
\hat V\;\; \text{is continuously differentiable on}\;\; [0,\infty)
^N\, .
\end{equation}
In particular, this shows that \eqref{7.3} holds as well as that
$\hat V$ belongs to $C^1(\bar C)$ as required in \eqref{7.4} above.
The fact that $b$ satisfies \eqref{7.5} was established in
\eqref{6.14}-\eqref{6.15} above. This shows that $(\hat V,b)$
belongs to $\cal C$ and solves \eqref{7.1}-\eqref{7.3} as claimed.
To derive uniqueness of the solution we will first see in the next
section that any solution $(U,a)$ to \eqref{7.1}-\eqref{7.3} from
the class $\cal C$ admits an explicit integral representation for
$U$ expressed in terms of $a$, which in turn solves a nonlinear
Fredholm integral equation, and we will see that this equation
cannot have other solutions satisfying the required properties. From
these facts we can conclude that the free-boundary problem
\eqref{7.1}-\eqref{7.3} cannot have other solutions in the class
$\cal C$ as claimed. This completes the proof. \hfill $\square$

\section{Nonlinear integral equations}

In this section we show that the optimal stopping boundary $b$ from
\eqref{6.16} can be characterised as the unique solution to a
nonlinear Fredholm integral equation. This also yields an explicit
integral representation of the value function $\hat V$ from
\eqref{3.4} expressed in terms of the optimal stopping boundary $b$.
As a consequence of the existence and uniqueness result for the the
nonlinear Fredholm integral equation we also obtain uniqueness of
the solution to the free-boundary problem \eqref{7.1}-\eqref{7.3} as
explained in the proof of Theorem 9 above. Finally, collecting the
results derived throughout the paper we conclude our exposition by
disclosing the solution to the initial problem.

\medskip

1.\ Let $p = p(t;\varphi_1, \ldots, \varphi_N,\psi_1, \ldots,
\psi_N)$ denote the transition probability density function of the
Markov process $\varPhi = (\varPhi^1, \ldots, \varPhi^N)$ in the
sense that
\begin{equation} \h{1pc} \label{8.1}
\PP_{\!\!\varphi_1, \ldots, \varphi_N}^\infty \big( \varPhi_t\!
\in\! A \big) = \int\! ... \! \int_A\, p(t;\varphi_1, \ldots,
\varphi_N, \psi_1, \ldots, \psi_N)\, d\psi_1 \ldots d\psi_N
\end{equation}
for any measurable $A \subseteq [0,\infty)^N$ with $(\varphi_1,
\ldots, \varphi_N) \in [0,\infty)^N$ and $t \ge 0$ given and fixed.
Having $p$ we can evaluate the expression of interest in the theorem
below as follows
\begin{align} \h{1pc} \label{8.2}
&K_b(t;\varphi_1, \ldots, \varphi_N) := \EE_{\varphi_1, \ldots,
\varphi_N}^\infty \big[ L(\varPhi_t^1, \ldots, \varPhi_t^N)\;\!
I \big(\varPhi_t^N\! <\! b(\varPhi_t^1, \ldots, \varPhi_t^{N-1})
\big) \big] \\[3pt] \notag &= \int\! ...\! \int_{\{\psi_N <
b(\psi_1, \ldots, \psi_{N-1})\}}\! L(\psi_1, \ldots, \psi_N)\,
p(t;\varphi_1, \ldots, \varphi_N,\psi_1, \ldots, \psi_N)\,
d\psi_1 \ldots d\psi_N
\end{align}
for $t \ge 0$ and $(\varphi_1, \ldots, \varphi_N)\in [0,\infty)^N$
where $L$ is defined in \eqref{4.3} above.

\medskip

\textbf{Theorem 10 (Existence and uniqueness).} \emph{The optimal
stopping boundary $b$ in \eqref{3.4} can be characterised as the
unique solution to the nonlinear Fredholm integral equation
\begin{equation} \h{4pc} \label{8.3}
\int_0^\infty\! e^{-\lambda t} K_b(t;\varphi_1, \ldots,
\varphi_{N-1},b(\varphi_1, \ldots, \varphi_{N-1}))\, dt = 0
\end{equation}
in the class of continuous {\rm \&} decreasing (convex) functions
$b$ on $[0,\infty)^{N-1}$ satisfying $\sum_{i=1}^{N-1} p_i$
$\varphi_i + p_N b(\varphi_1, \ldots, \varphi_{N-1}) \m \lambda/c
\ge 0$ for $(\varphi_1, \ldots, \varphi_{N-1}) \in
[0,\infty)^{N-1}$. The value function $\hat V$ in \eqref{3.4} admits
the following representation
\begin{equation} \h{4pc} \label{8.4}
\hat V(\varphi_1, \ldots, \varphi_N) = \int_0^\infty\! e^{-\lambda
t} K_b(t; \varphi_1, \ldots, \varphi_N)\, dt
\end{equation}
for $(\varphi_1, \ldots, \varphi_N)\in [0,\infty)^N$. The optimal
stopping time in \eqref{3.4} is given by
\begin{equation} \h{4pc} \label{8.5}
\tau_b = \inf\, \{\, t \ge 0\; \vert\; \varPhi_t^N \ge b(\varPhi
_t^1, \ldots, \varPhi_t^{N-1})\, \}
\end{equation}
under $\PP_{\!\!\varphi_1, \ldots, \varphi_N}^\infty$ with
$(\varphi_1, \ldots, \varphi_N)\in [0,\infty)^N$ given and fixed.}

\medskip

\textbf{Proof.} (i) \emph{Existence}. We first show that the optimal
stopping boundary $b$ in \eqref{3.4} solves \eqref{8.3}. Recalling
that $b$ satisfies the properties stated following \eqref{6.14}
above, this will establish the existence of a solution to
\eqref{8.3} in the specified class of functions.

For this, to gain control over the (individual) second partial
derivatives $\hat V_{\varphi_i \varphi_j}$ close to the optimal
stopping boundary within $C$ for $1 \le i,j \le N$ (see \cite{GT}
for general results of this kind in the elliptic case), consider the
sets $C_n := \{\, \varphi \in [0,\infty)^N\; \vert\; \hat V(\varphi)
< -1/n\, \}$ and $D_n := \{\, \varphi \in [0,\infty)^N\; \vert\;
\hat V(\varphi) \ge -1/n\, \}$ for $n \ge 1$ (large). Note that $C_n
\uparrow C$ and $D_n \downarrow D$ as $n \uparrow \infty$. Moreover,
using the same arguments as for the sets $C$ and $D$ above, we find
that the set $D_n$ is convex, and the boundary $b_n = b_n(\varphi_1,
\ldots, \varphi_{N-1})$ between $C_n$ and $D_n$ is a convex,
continuous, decreasing function of $(\varphi_1, \ldots,
\varphi_{N-1})$ in $[0,\infty)^{N-1}$. This also shows that $b_n
\uparrow b$ uniformly on $[0,\infty)^{N-1}$ as $n \rightarrow
\infty$ (where the functions $b_n$ and $b$ take zero value at
$(\varphi_1, \ldots, \varphi_{N-1})$ by definition whenever
$(\varphi_1, \ldots, \varphi_{N-1},0)$ belongs to $D_n$ and $D$
respectively for $n \ge 1$). Approximate the value function $\hat V$
in \eqref{3.4} by functions $\hat V^n$ defined as $\hat V$ on $C_n$
and $-1/n$ on $D_n$ for $n \ge 1$. Note that $\hat V^n \uparrow \hat
V$ uniformly on $[0,\infty)^N$ as $n \rightarrow \infty$. Moreover,
letting $n \ge 1$ be given and fixed in the sequel, clearly $\hat
V^n$ is a continuous function on $[0,\infty)^N$ and $\hat V^n$
restricted to $C_n$ and $D_n$ belongs to $C^2(\bar C_n)$ and
$C^2(\bar D_n)$ respectively. Finally, since $b_n$ is convex we know
that $b_n(\varPhi^1, \ldots, \varPhi^{N-1})$ is a continuous
semimartingale. This shows that the change-of-variable formula with
local time on surfaces \cite[Theorem 2.1]{Pe-3} is applicable to
$\hat V^n$ composed with $\varPhi = (\varPhi^1, \ldots, \varPhi^N)$
and using \eqref{7.1} this gives
\begin{align} \h{-0.5pc} \label{8.6}
e^{-\lambda t}\:\! \hat V^n(\varPhi_t) &= \hat V^n(\varPhi_0)
+ \int _0^t e^{-\lambda s} \big( \LL_\varPhi \hat V^n \m \lambda
\hat V^n \big)(\varPhi_s)\, ds \\[-3pt] \notag &\h{13pt}+ \sum_{i=1}^N
\sum_{j=1}^k \int_0^t e^{-\lambda s}\:\! \hat V^n_{\varphi_i}
(\varPhi_s)\:\! \mu\;\! \varPhi_s^i\: dB_s^{n_j(i)} - \int_0^t
e^{-\lambda s}\:\! \hat V^n _{\varphi_N}(\varPhi_s-)\, d\ell_s^{b_n}
(\varPhi) \\ \notag &= \hat V^n(\varPhi_0) - \int_0^t e^{-\lambda
s}\:\! L(\varPhi_s)\;\! I(\varPhi_s\! \in\! C_n)\, ds + M_t^n -
\int_0^t e^{-\lambda s}\:\! \hat V_{\varphi_N}(\varPhi_s)\,
d\ell_s^{b_n}(\varPhi)
\end{align}
where $M_t^n = \sum_{i=1}^N \sum_{j=1}^k \int_0^t e^{-\lambda s}\:\!
\hat V_{\varphi_i} (\varPhi_s)\:\! \mu\:\! \varPhi_s^i\;\!
I(\varPhi_s\! \in\! C_n)\: dB_s^{n_j(i)}$ is a continuous martingale
for $t \ge 0$ and $\ell^{b_n}(\varPhi)$ is the local time of
$\varPhi$ on the curve $b_n$ given by
\begin{align} \h{2pc} \label{8.7}
\ell_t^{b_n}(\varPhi) = \PP \text{-}\! \lim_{\eps \downarrow
0} \frac {1}{2 \eps}\! \int_0^t &I( -\eps < \varPhi_s^N\! \m
b_n(\varPhi_s^1, \ldots, \varPhi_s^{N-1}) < \eps ) \\[-2pt] \notag
&d \langle \varPhi^N\! \m b_n(\varPhi^1, \ldots, \varPhi^{N-1}),
\varPhi^N\! \m b_n(\varPhi^1, \ldots, \varPhi^{N-1}) \rangle_s
\end{align}
for $t \ge 0$. To gain control over the final term in \eqref{8.6},
note that the It\^o-Tanaka formula (cf.\ \cite[pp 222-223]{RY})
yields
\begin{align} \h{-0.5pc} \label{8.8}
&\big( b_n(\varPhi_t^1, \ldots, \varPhi_t^{N-1}) \m \varPhi_t^N
\big)^+ = \big( b_n(\varPhi_0^1, \ldots, \varPhi_0^{N-1}) \m
\varPhi_0^N \big)^+ \\ \notag &\h{13pt}+\! \int_0^t I(b_n(
\varPhi_s^1, \ldots, \varPhi_s^{N-1}) \m \varPhi_s^N\! >\! 0)
\, d( b_n(\varPhi^1, \ldots, \varPhi^{N-1}) \m \varPhi^N)_s +
\frac{1}{2}\;\! \ell_t^{b_n} (\varPhi) \\ \notag &= \big(
b_n(\varPhi_0^1, \ldots, \varPhi_0^{N-1}) \m \varPhi_0^N \big)^+
\! +\! \int_0^t I(b_n(\varPhi_s^1, \ldots, \varPhi_s^{N-1}) \m
\varPhi_s^N\! >\! 0) \\[-2pt] \notag &\h{16pt} \Big( \sum_{i=1}
^{N-1} \frac{\partial b_n}{\partial \varphi_i}(\varPhi_s^1,
\ldots \varPhi_s^{N-1})\;\! d \varPhi_s^i + \frac{1}{2} \sum
_{i=1}^{N-1} \sum_{j=1}^{N-1} \frac{\partial^2 b_n}{\partial
\varphi_i \partial \varphi_j}(\varPhi_s^1, \ldots \varPhi_s
^{N-1})\, d \langle \varPhi^i,\varPhi^j \rangle_s - d \varPhi
_s^N \Big) \\ \notag &\h{13pt}+ \frac{1}{2}\;\! \ell_t^{b_n}
(\varPhi)
\end{align}
for $t \ge 0$ where we use that $b_n$ is $C^2$ by the implicit
function theorem since the smooth fit fails at $b_n$ due to its
suboptimality in the problem \eqref{3.4}. Since $b_n$ is convex we
know that the Hessian matrix of $b_n$ is non-negative definite and
recalling from \eqref{3.3}+\eqref{4.1} that $d \langle
\varPhi^i,\varPhi^j \rangle_s = \mu^2 \varPhi_s^i \varPhi_s^j
(I_i,I_j)\, ds$ for $1 \le i,j \le N-1$ we see that the integral
associated with the double sum in \eqref{8.8} is non-negative. It
follows therefore from \eqref{8.8} using \eqref{3.3} above that
\begin{align} \h{-0.5pc} \label{8.9}
\frac{1}{2}\;\! \ell_t^{b_n}(\varPhi) &\le \big( b_n(\varPhi_t^1,
\ldots, \varPhi_t^{N-1}) \m \varPhi_t^N \big)^+ \\ \notag &\h{14pt}
- \sum_{i=1}^{N-1} \int_0^t I(b_n(\varPhi_s^1, \ldots, \varPhi_s
^{N-1}) \m \varPhi_s^N\! >\! 0)\, \frac{\partial b_n}{\partial
\varphi_i} (\varPhi_s^1, \ldots, \varPhi_s^{N-1}) \;\! \lambda
(1 \p \varPhi_s^i) \, ds \\ \notag &\h{14pt} +\! \int_0^t
I(b_n(\varPhi_s^1, \ldots, \varPhi_s^{N-1}) \m \varPhi_s^N\!
>\! 0)\;\! \lambda (1 \p \varPhi_s^N) \, ds + N_t^n
\end{align}
where \v{-2pt} $N_t^n = - \sum_{i=1}^{N-1} \sum_{j=1}^k \int_0^t
I(b_n(\varPhi_s^1, \ldots, \varPhi_s ^{N-1}) - \varPhi_s^N\! >\! 0)
(\partial b_n/\partial \varphi_i)(\varPhi_s^1, \ldots, \varPhi_s
^{N-1})\;\! \mu\;\! \varPhi_s^i\:\! dB_s^{n_j(i)}$ $+ \sum_{j=1}^k
\int_0^t I(b_n(\varPhi_s^1, \ldots, \varPhi_s ^{N-1})-\varPhi_s^N \!
>\! 0)\;\! \mu\;\! \varPhi_s^N\:\! dB_s^{n_j(N)}$ is a continuous
local martingale for $t \ge 0$. Let $(\tau_m)_{m \ge 1}$ be a
localising sequence of stopping times for $N^n$, define the stopping
time
\begin{equation} \h{5pc} \label{8.10}
\sigma_m = \inf\, \{\, t \ge 0\; \vert\; (\varPhi_t^1, \ldots,
\varPhi_t^{N-1}) \notin (\tfrac{1}{m},\infty)^{N-1}\, \}
\end{equation}
and set $\rho_m := \tau_m \wedge \sigma_m$ for $m \ge 1$. Fixing
$\varphi \in [0,\infty)^N$ we then find from \eqref{8.9} that
\begin{align} \h{3pc} \label{8.11}
\frac{1}{2}\;\! \EE_\varphi^\infty \big[ \ell_{t \wedge \rho_m}
^{b_n} (\varPhi) \big] &\le \varphi_N^* - \sum_{i=1}^{N-1} \frac
{\partial b_n}{\partial \varphi_i}(\tfrac{1}{m}, \ldots, \tfrac{1}
{m}) \int_0^t \lambda \big( 1 \p \EE_\varphi^\infty (\varPhi_s^i)
\big)\, ds \\ \notag &\h{13pt}+ \int_0^t \lambda \big( 1 \p
\EE_\varphi^\infty (\varPhi_s^N) \big)\, ds \le K_m(t)
\end{align}
for $t \ge 0$ and $m \ge 1$ where the positive constant $K_m(t)$
does not depend on $n \ge 1$ because each $b_n$ is convex and $b_n
\uparrow b$ on $[0,1/m]^{N-1}$ as $n \rightarrow \infty$ so that
$(\partial b_n/\partial \varphi_i)(1/m, \ldots, 1/m)$ must stay
bounded from below over $n \ge 1$ for $1 \le i \le N \m 1$ if $b_n$
is to stay below $b$ on $[0,1/m]^{N-1}$ for all $n \ge 1$. In
addition, by \eqref{7.7} we know that $\hat V_{\varphi_N}$ is
continuous on $\bar C$ and hence uniformly continuous too because
$\bar C$ is a compact set. \v{-1pt} It follows therefore that $0 \le
\hat V_{\varphi_N}(\varphi_1, \ldots, \varphi_{N-1},b_n(\varphi_1,
\ldots, \varphi_{N-1})) \le \eps$ for all $(\varphi_1, \ldots,
\varphi_{N-1}) \in [0,\infty)^{N-1}$ and all $n \ge n_\eps$ with
$n_\eps \ge 1$ large enough depending on the given and fixed
$\eps>0$. Combining this fact with \eqref{8.11}, upon replacing $t$
with $t \wedge \rho_m$ in the final integral of \eqref{8.6} and
taking $\EE_\varphi^\infty$ of the resulting expression, we see that
\begin{equation} \h{5pc} \label{8.12}
0 \le \EE_\varphi^\infty \bigg[ \int_0^{t \wedge \rho_m}\!\!
e^{-\lambda s}\:\! \hat V_{\varphi_N}(\varPhi_s)\, d\ell_s^
{b_n}(\varPhi) \bigg] \le 2\:\! \eps\:\! K_t(m)
\end{equation}
for all $n \ge n_\eps$ with $t \ge 0$ and $m \ge 1$ given and fixed.
This shows that the expectation in \eqref{8.12} tends to zero as $n$
tends to infinity for every $t \ge 0$ and $m \ge 1$ given and fixed.
Using this fact in \eqref{8.6} upon replacing $t$ with $t \wedge
\rho_m$, taking $\EE_\varphi^\infty$ on both sides, and letting $n$
tend to infinity, we find by the monotone convergence theorem upon
recalling \eqref{7.6} that
\begin{equation} \h{3pc} \label{8.13}
\hat V(\varphi) = \EE_\varphi^\infty \big[ e^{-\lambda(t \wedge
\rho _m)} \:\! \hat V(\varPhi_{t \wedge \rho_m}) \big] + \EE_
\varphi^\infty \bigg[\! \int_0^{t \wedge \rho_m}\!\! e^{-\lambda
s}\:\! L(\varPhi_s) \:\! I(\varPhi_s\! \in\! C) \, ds \bigg]
\end{equation}
for all $t \ge 0$ and all $m \ge 1$. Letting $m \rightarrow \infty$
and using that $\rho_m \rightarrow \infty$ because $0$ is an
entrance boundary point for each $\varPhi^i$ with $1 \le i \le N \m
1$, we see from \eqref{8.13} upon recalling \eqref{7.6} and using
the dominated convergence theorem that
\begin{equation} \h{3pc} \label{8.14}
\hat V(\varphi) = \EE_\varphi^\infty \big[ e^{-\lambda t}\:\!
\hat V(\varPhi_t) \big] + \EE_\varphi^\infty \bigg[\! \int_0^t
e^{-\lambda s} \:\! L(\varPhi _s)\:\! I(\varPhi_s\! \in\! C)
\, ds \bigg]
\end{equation}
for all $t \ge 0$. Finally, letting $t \rightarrow \infty$ in
\eqref{8.14} and using the dominated and monotone convergence
theorems upon recalling \eqref{7.6}, we find that
\begin{equation} \h{5pc} \label{8.15}
\hat V(\varphi) = \EE_\varphi^\infty \bigg[\! \int_0^\infty
e^{-\lambda s}\:\! L(\varPhi _s)\:\! I(\varPhi_s\! \in\! C)
\, ds
\bigg]
\end{equation}
for all $\varphi \in [0,\infty)^N$. Recalling \eqref{6.14} and
\eqref{8.2} above we see that this establishes the representation
\eqref{8.4} as claimed. Moreover, the fact that $\tau_b$ from
\eqref{8.5} is optimal in \eqref{3.4} follows by \eqref{6.14} above.
Finally, inserting $\varphi_N = b(\varphi_1, \ldots, \varphi_{N-1})$
in \eqref{8.4} and using that $\hat V(\varphi_1, \ldots,
\varphi_{N-1},b(\varphi_1, \ldots, \varphi_{N-1}))=0$, we see that
$b$ solves \eqref{8.3} as claimed.

\smallskip

(ii) \emph{Uniqueness}. To show that $b$ is a unique solution to the
equation \eqref{8.3} in the specified class of functions, one can
adopt the four-step procedure from the proof of uniqueness given in
\cite[Theorem 4.1]{DuPe} extending and further refining the original
uniqueness arguments from \cite[Theorem 3.1]{Pe-1}. Given that the
present setting creates no additional difficulties we will omit
further details of this verification and this completes the proof.
\hfill $\square$

\medskip

2.\ The nonlinear Fredholm integral equation \eqref{8.3} can be used
to find the optimal stopping boundary $b$ numerically (using Picard
iteration). Inserting this $b$ into \eqref{8.4} via \eqref{8.2} we
also obtain a closed form expression for the value function $\hat
V$. Collecting the results derived throughout the paper we now
disclose the solution to the initial problem.

\medskip

\textbf{Corollary 11.} \emph{The value function in the initial problem
\eqref{2.4} is given by
\begin{equation} \h{5pc} \label{8.16}
V(\pi) = (1 \m \pi)\, \Big[\:\! 1 + c\;\! \hat V \Big( \frac{\pi}{1
\m \pi}, \ldots, \frac{\pi}{1 \m \pi} \Big)\:\! \Big]
\end{equation}
for $\pi \in [0,1]$ where the function $\hat V$ is given by
\eqref{8.4} above. The optimal stopping time in the initial problem
\eqref{2.4} is given by
\begin{align} \h{-1pc} \label{8.17}
\tau_* = \inf \Big\{ t \ge 0\: &\big \vert\: e^{\;\! \mu \sum_{j=1}^k
\! X_t^{n_j(N)} \! + (\lambda - k \frac{\mu^2}{2})\:\! t} \Big(
\tfrac{\pi}{1-\pi} \p \lambda \!\! \int_0^t\! e^{\;\! - \mu
\sum_{j=1}^k\! X_s^{n_j(N)}\! - (\lambda - k \frac{\mu^2} {2})\:\!
s}\, ds \Big) \\ \notag &\h{6pt}\ge b \Big( e^{\;\! \mu
\sum_{j=1}^k\! X_t^{n_j(1)} \! + (\lambda - k \frac{\mu^2}{2})\:\!
t} \Big( \tfrac{\pi}{1-\pi} \p \lambda \!\! \int_0^t\! e^{\;\! -\mu
\sum_{j=1}^k\! X_s^{n_j(1)}\! - (\lambda - k \frac{\mu^2} {2})\:\!
s}\, ds \Big), \ldots, \\ \notag &\h{36pt}e^{\;\! \mu \sum_{j=1}^k\!
X_t^{n_j(N-1)} \! + (\lambda - k \frac{\mu^2}{2})\:\! t} \Big(
\tfrac{\pi}{1-\pi} \p \lambda \!\! \int_0^t\! e^{\;\! -\mu
\sum_{j=1}^k\! X_s^{n_j(N-1)}\! - (\lambda - k \frac{\mu^2} {2})\:\!
s}\, ds \Big) \Big) \Big\}
\end{align}
where $b$ is a unique solution to \eqref{8.3} above (see Figure 1).}

\textbf{Proof.} The identity \eqref{8.16} was established in
\eqref{3.1} above. The explicit form of the optimal stopping time
\eqref{8.17} follows from \eqref{8.5} in Theorem 10 combined with
\eqref{2.11}+\eqref{2.12} above. The final claim on $b$ was derived
in Theorem 10 above. This completes the proof. \hfill $\square$

\section{General case}

In the general case of the quickest detection problem \eqref{2.4} we
no longer insist that exactly $k$ of the coordinate processes $X^1,
\ldots, X^n$ get a (known) non-zero drift $\mu$ but instead allow
that \emph{any} number of them get such a drift with prescribed
probabilities. In this section we show that the methodology
developed in the previous sections to solve the problem for exactly
$k$ coordinate processes can be used to solve the problem in the
general case for \emph{any} number of coordinate processes. This
extension of the solution will also enable us to reduce the
dimension of the problem to its minimal value which is of
fundamental importance in real applications.

1.\ The key issue in extending the solution from exactly $k$ to any
number of the coordinate processes is whether the hypoelliptic
structure established in Section 5 above using the H\"ormander
theorem remains preserved. For this, we first note that this is the
case if we mechanically increase the problem dimension to a size
which however could be alarmingly high. For example, returning to
the case when $n=10$ and $k=5$ as discussed in Section 1 above, and
allowing any (one or several) tagged $4$ coordinate processes to get
the drift as well, we would increase the problem dimension from $N_1
= \binom{10}{5} = 252$ to $N_2 = \binom{10}{5} + \binom{10}{4} =
462$. Clearly, the larger the number of the tagged 4 coordinate
processes involved, the more justified increase of the problem
dimension would be, and vice versa. This raises the question of
establishing a minimal dimension of the problem given all the tagged
coordinate processes that can get such a drift with prescribed
probabilities. On closer inspection of the previous arguments we
then note that the H\"ormander condition, and hence the hypoelliptic
structure established in Section 5 above as well, remain valid in
the case of a minimal dimension of the problem. This makes the
results derived in the previous sections applicable in the general
case where the number of the tagged coordinate processes is firstly
enlarged to an arbitrary value and then trimmed down to the value of
a minimal problem dimension.

\medskip

2.\ To describe the solution to the quickest detection problem
\eqref{2.4} in the general case with a minimal problem dimension, we
will return to the beginning of Section 2 and replace the number $1
\le k \le n$ by the numbers $1 \le k_1 < \dots < k_m \le n$ with $m
\le n$. This means that any $k_l$ of the coordinate processes $X^1,
\ldots, X^n$ get a (known) non-zero drift $\mu$ at time $\theta$ for
$1 \le l \le m$ instead of exactly $k$ of them. Setting $C_{k_1,
\ldots, k_m}^{\:\!n} := \cup_{l=1}^m C_{k_l}^{\:\!n}$ we see that
that the random variable $\beta$ in \eqref{2.1} taking values in the
set $C_{k_1, \ldots, k_m}^{\:\!n}$ satisfies $\PP_{\!\pi}(\beta =
(n_1, \ldots, n_{k_l})) = p_{n_1, \ldots, n_{k_l}}$ \v{-1pt} for
some $p_{n_1, \ldots, n_{k_l}} \in [0,1]$ with $\sum_{l=1}^m \sum
p_{n_1, \ldots, n_{k_l}} = 1$ given and fixed where the second sum
is taken over all $(n_1, \ldots, n_{k_l}) \in C_{k_l}^{\:\!n}$. It
is important in this setting that at least one among $p_{n_1,
\ldots, n_{k_l}}$ when $(n_1, \ldots, n_{k_l})$ runs through
$C_{k_l}^{\:\!n}$ is assumed to be strictly positive for every $
1\le l \le m$ given and fixed. If this would not be the case for
some $1 \le l \le m$ then $C_{k_l}$ could be omitted from the
setting. As before, with a slight abuse of notation, in \eqref{2.1}
we write $i \in \beta$ to express the fact that $i$ belongs to the
set $\{ n_1, \ldots, n_{k_l} \}$ consisting of the elements which
form $\beta = (n_1, \ldots, n_{k_l})$ in $C_{k_l}^{\:\!n}$ for $1
\le l \le m$. This means that $n_1, \ldots, n_{k_l} \in \beta$ if
and only if the coordinate processes $X^{n_1}, \ldots, X^{n_{k_l}}$
get drift $\mu$ at time $\theta$ with probability $p_{n_1, \ldots,
n_{k_l}}$ for $(n_1, \ldots, n_{k_l}) \in C_{k_l}^{\:\!n}$ with $1
\le l \le m$. With a similar abuse of notation, which will be
helpful in what follows as before, we will first arrange the
elements of $C_{k_1, \ldots, k_m}^{\:\!n}$ in a lexicographic order
starting first with $C_{k_1}^{\:\!n}$ and moving forward until we
reach $C_{k_m}^{\:\!n}$, then remove all $(n_1, \ldots, n_{k_l})$
from the ordered set $C_{k_1, \ldots, k_m}^{\:\!n}$ for which
$p_{n_1, \ldots, n_{k_l}} = 0$ when $1 \le l \le m$ and denote
\v{-3pt} the remaining (ordered) set by $\hat C_{k_1, \ldots,
k_m}^{\:\!n}$, and finally identify the $i$\!-th element of the
ordered set $\hat C_{k_1, \ldots, k_m}^{\:\!n}$ by its index $i$
itself for $1 \le i \le N$ where we let $N$ denote the total number
of elements in $\hat C_{k_1, \ldots, k_m}^{\:\!n}$. \v{-3pt} Thus as
before we write $i = (n_1, \ldots, n_{k_l}) \in \hat C_{k_1, \ldots,
k_m}^{\:\!n}$ for $1 \le l \le m$ to express \h{-5pt} this
identification explicitly for $1 \le i \le N$. Note that
\begin{equation} \h{3pc} \label{9.1}
\hat C_{k_1, \ldots, k_m}^{\:\!n} = \{\, (n_1, \ldots, n_{k_l})
\in C_{k_l}^{\:\!n} \; \vert\; p_{n_1, \ldots, n_{k_l}}\! > 0
\;\; \text{for} \;\; 1 \le l \le m\, \}
\end{equation}
and \v{-1pt} hence $N \le \binom{n}{k_1} + \ldots + \binom{n}{k_m}$
with equality being attained if $p_{n_1, \ldots, n_{k_l}} > 0$ for
all $(n_1, \ldots, n_{k_l})$ $\in C_{k_1, \ldots, k_m}^{\:\!n}$\!
with $1 \le l \le m$.

\medskip

3.\ The rest of the analysis in Sections 2-4 above can then be
carried out in exactly the same way. The only difference is
notational in that the index set $C_k^n$ needs to be firstly
enlarged and then trimmed down to the index set $\hat C_{k_1,
\ldots, k_m}^{\:\!n}$ as explained above. In particular, from the
enlarged/trimmed identities \eqref{2.11}+\eqref{2.12} using It\^o's
formula we find that the enlarged/trimmed system of stochastic
differential equations \eqref{3.3} reads as follows
\begin{equation} \h{0.5pc} \label{9.2}
d \varPhi_t^i = \lambda\:\! (1 \p \varPhi_t^i)\;\! dt +
\sum_{j=1}^{k_l} \mu\;\! \varPhi_t^i\, dB_t^{n_j}\;\;
(i\! =\! (n_1, \ldots, n_{k_l})\! \in\! \hat C_{k_1,
\ldots, k_m}^{\:\!n}\; \text{with}\;\, 1 \le l \le m)
\end{equation}
under $\PP^\infty$ with $\varPhi_0^i = \varphi_i$ in $[0,\infty)$
all being equal to $\pi/(1 \m \pi)$ for $1 \le i \le N$ in $\hat
C_{k_1, \ldots, k_m}^{\:\!n}$ and $\pi \in [0,1)$. Similarly, we see
that the enlarged/trimmed optimal stopping problem \eqref{3.4}
remains unchanged and reads
\begin{equation} \h{5pc} \label{9.3}
\hat V(\varphi) = \inf_\tau\;\! \EE_\varphi^\infty \Big[ \int
_0^\tau e^{-\lambda t} \Big( \sum_{i=1}^N p_i\;\! \varPhi_t^i
- \frac{\lambda} {c}\;\! \Big)\;\! dt\:\! \Big]
\end{equation}
for $\varphi \in [0,\infty)^N$ with $\PP_{\!\varphi}^\infty
(\varPhi_0\! =\! \varphi)=1$ \v{-1pt} where the infimum is taken
over all stopping times $\tau$ of $\varPhi$ and we recall that $p_i
\in (0,1]$ for $1 \le i \le N$ in $\hat C_{k_1, \ldots,
k_m}^{\:\!n}$ with $\sum_{i=1}^N p_i = 1$.

\medskip

4.\ This raises the question whether the hypoelliptic structure
established using the H\"orman- der theorem in Section 5 above is
preserved in the enlarged/trimmed system \eqref{9.2} above. We now
show that the answer is affirmative.

\medskip

\textbf{Corollary 12.} \emph{The results of Proposition 3 and
Proposition 4 above remain valid for the enlarged/trimmed system
\eqref{9.2} above.}

\smallskip

\textbf{Proof.} We show that the H\"ormander condition \eqref{5.8}
is satisfied using backward induction over $k_1 < \ldots < k_m$.
Setting $k=k_m$ and connecting to the conclusion \eqref{5.17} in the
proof of Proposition 3 above, we see that $\partial_{\varphi_j} \in
Lie\:\! (D_0,D_1, \ldots, D_N)$ for $N_{m-1} \p 1 \le j \le N_m$
where $N_l$ denotes the number of elements from $\hat C_{k_1,
\ldots, k_m}^{\:\!n}$ that belong to $C_{k_l}^{\:\!n}$ for $1 \le l
\le m$. Note that for that conclusion we did not use that $N_l$ must
be equal to its general bound $\binom{n}{k_l}$ but could also be any
strictly smaller number as well for $1 \le l \le m$. Recall also
that the lexicographic order of $\hat C_{k_1, \ldots, k_m}^{\:\!n}$
is without loss of generality assumed to start with
$C_{k_1}^{\:\!n}$ and end with $C_{k_m}^{\:\!n}$. The induction step
can therefore be realised by setting $k = k_{m-1}$. Then the same
arguments as in \eqref{5.13}-\eqref{5.15} show that
\begin{equation} \h{3pc} \label{9.4}
[\:\![\:\![D_0,D_{n_1}],D_{n_2}], \ldots, D_{n_{k_{m-1}}}] \sim
\sum_{i=1}^N I_{i\:\!n_1} I_{i\:\!n_2} \ldots I_{i\:\!n_{k_{m-1}}}
\:\! \partial_{\varphi_i}
\end{equation}
for $j=(n_1, \ldots, n_{k_{m-1}}) \in \hat C_{k_1, \ldots,
k_m}^{\:\!n}$ given and fixed. This reveals a crucial difference in
comparison with \eqref{5.16} above because this time we have
\begin{align} \h{5pc} \label{9.5}
I_{i\:\!n_1} I_{i\:\!n_2} \ldots I_{i\:\!n_{k_{m-1}}} &= 1\;\;
\text{if}\;\; i = j\;\; \text{or}\;\; i \in \hat C_{k_m}^{\:\!n}(j)
\\[3pt] \notag &= 0\;\; \text{otherwise}
\end{align}
where $\hat C_{k_m}^{\:\!n}(j)$ denotes the set of all $i \in
C_{k_m}^{\:\!n}$ in $\hat C_{k_1, \ldots, k_m}^{\:\!n}$ which
include $j$ in the sense that $n_1, \ldots, n_{k_{m-1}}$ belong to
$i$ understood as the set of its elements. It follows therefore from
\eqref{9.4} and \eqref{9.5} that we have
\begin{equation} \h{5pc} \label{9.6}
[\:\![\:\![D_0,D_{n_1}],D_{n_2}], \ldots, D_{n_{k_{m-1}}}] \sim
\partial_{\varphi_j} +\!\! \sum_{i \in \hat C_{k_m}^{\:\!n}(j)}
\partial_{\varphi_i}\, .
\end{equation}
Recalling however that $\partial_{\varphi_i} \in Lie\:\! (D_0,D_1,
\ldots, D_N)$ for every $i \in \hat C_{k_m}^{\:\!n}(j)$ as
established in the previous step, we see from \eqref{9.6} that
$\partial_{\varphi_j} \in Lie\:\! (D_0,D_1, \ldots, D_N)$ for
$N_{m-2} \p 1 \le j \le N_{m-1}$. Continuing this procedure by
backward induction until reaching $k = k_1$ we find that
$\partial_{\varphi_j} \in Lie\:\! (D_0,D_1, \ldots, D_N)$ for all $1
\le j \le N$ so that the H\"ormander condition \eqref{5.8} is
satisfied as claimed. The remaining arguments are identical to those
presented in the proofs of Proposition 3 and Proposition 4 above.
This completes the proof. \hfill $\square$

\medskip

5.\ The rest of the analysis in Sections 6-8 above can then be
carried out in exactly the same way with the index set $C_k^{\:\!n}$
replaced by the enlarged/trimmed index set $\hat C_{k_1, \ldots,
k_m}^{\:\!n}$. This yields the solution to the initial problem when
the number of coordinate processes getting a (known) non-zero drift
$\mu$ with prescribed (non-zero) probabilities is no longer fixed to
be $k$ exactly but could be \emph{any} number instead.

\medskip

\textbf{Corollary 13.} \emph{The results of Theorem 10 and Corollary
11 remain valid in the general case described by
\eqref{9.1}-\eqref{9.3} above.}

\medskip

\textbf{Proof.} Having established the hypoelliptic structure of
\eqref{9.1}-\eqref{9.3} in Corollary 12 above we can carry out the
proofs of Theorem 10 and Corollary 11 in exactly the same way. Note
that this also applies to all other results in Sections 6-8 (as well
as to those in Sections 2-4). This completes the proof. \hfill
$\square$

\medskip

\textbf{Acknowledgements.} The authors gratefully acknowledge
support from the United States Army Research Office Grant
ARO-YIP-71636-MA. The third-named author is indebted to Xinfu Chen
and David Saunders for the insightful discussions on hypoellipticity
in relation to \cite{JP} and related equations during a splendid
stay in Beijing, China (June/July 2017).

\begin{center}

\end{center}


\begin{thebibliography}{99}

\setlength{\itemsep}{2pt} \vspace{4pt}

\bibitem{BG} \textsc{Blumenthal, R. M. \emph{and} Getoor, R. K.}
(1968). \emph{Markov Processes and Potential Theory}. Academic
Press.

\bibitem{DePe} \textsc{De Angelis, T. \emph{and} Peskir, G.} (2020).
Global $C^1$ regularity of the value function in optimal stopping
problems. \emph{Ann. Appl. Probab.} 30 (1007--1031).

\bibitem{DuPe} \textsc{Du Toit, J. \emph{and} Peskir, G.} (2009).
Selling a stock at the ultimate maximum. \emph{Ann. Appl. Probab.}
19 (983--1014).

\bibitem{EP} \textsc{Ernst, P. A. \emph{and} Peskir, G.} (2022).
Quickest real-time detection of a Brownian coordinate drift.
\emph{Ann. Appl. Probab.} 32 (2652--2670).

\bibitem{GT} \textsc{Gilbarg, D. \emph{and} Trudinger, N. S.}
(2001). \emph{Elliptic Partial Differential Equations of the Second
Order}. Springer.

\bibitem{Ho} \textsc{H\"ormander, L.} (1967). Hypoelliptic second
order differential equations. \emph{Acta Math.} 119 (147--171).

\bibitem{JP} \textsc{Johnson, P. \emph{and} Peskir, G.} (2017).
Quickest detection problems for Bessel processes. \emph{Ann. Appl.
Probab.} 27 (1003--1056).

\bibitem{Pe-1} \textsc{Peskir, G.} (2005). On the American option
problem. \emph{Math. Finance} 15 (169--181).

\bibitem{Pe-2} \textsc{Peskir, G.} (2006). On the fundamental
solution of the Kolmogorov-Shiryaev equation. \emph{The Shiryaev
Festschrift} (Metabief, 2005), Springer (535--546).

\bibitem{Pe-3} \textsc{Peskir, G.} (2007). A change-of-variable
formula with local time on surfaces. \emph{S\'em. de Probab.} XL,
\emph{Lecture Notes in Math.} 1899, Springer (69--96).

\bibitem{Pe-4} \textsc{Peskir, G.} (2022). Weak solutions in the
sense of Schwartz to Dynkin's characteristic operator equation.
\emph{Research Report} No. 1, \emph{Probab. Statist. Group
Manchester} (20 pp).

\bibitem{PS} \textsc{Peskir, G. \emph{and} Shiryaev, A. N.} (2006).
\emph{Optimal Stopping and Free-Boundary Problems}. Lectures in
Mathematics, ETH Z\"{u}rich, Birkh\"{a}user.

\bibitem{RY} \textsc{Revuz, D. \emph{and} Yor, M.} (1999).
\emph{Continuous Martingales and Brownian Motion.} Springer-Verlag.

\bibitem{RW} \textsc{Rogers, L. C. G. \emph{and} Williams, D.}
(2000). \emph{Diffusions, Markov Processes and Martingales: It\^o
Calculus} (Vol 2). Cambridge University Press.

\bibitem{Sh-1} \textsc{Shiryaev, A. N.} (1961). The problem of the
most rapid detection of a disturbance in a stationary process.
\emph{Soviet Math. Dokl.} 2 (795--799).

\bibitem{Sh-2} \textsc{Shiryaev, A. N.} (1978). \emph{Optimal
Stopping Rules.} Springer-Verlag.

\bibitem{Sh-3} \textsc{Shiryaev, A. N.} (1996). \emph{Probability}.
Springer.

\bibitem{Sh-4} \textsc{Shiryaev, A. N.} (2010). Quickest detection
problems: Fifty years later. \emph{Sequential Anal.} 29 (345--385).

\end{thebibliography}
\end{document}